\documentclass[a4paper,10pt]{amsart}
\usepackage{epsf}  
\usepackage{amsfonts, amssymb, amsthm, amsmath}  
\usepackage{hyperref}
\usepackage{graphicx}
\usepackage{psfrag}
\newtheorem{thm}{Theorem}  

\newtheorem{cor}[thm]{Corollary}  
\newtheorem{lemma}[thm]{Lemma}  
\newtheorem{remark}[thm]{Remark}  
\newtheorem{defn}[thm]{Definition}

\numberwithin{thm}{section}  
\def\pf{\noindent\emph{Proof: }}  
\def\stop{\hfill$\square$}

\providecommand{\totl}[1]{\ensuremath{\lceil #1\rceil }}
\providecommand{\totb}[1]{\ensuremath{\underline{ #1}}}
\DeclareMathOperator{\End}{End}
\DeclareMathOperator{\Aut}{Aut}

\newcommand{\ex}{\bold}
\providecommand {\e}[1]{\mathfrak t^{#1}}

\DeclareMathOperator{\expl}{Expl}

\DeclareMathOperator{\trop}{Trop}
\newcommand{\dbar}{\bar{\partial}}

\providecommand{\et}[2]{\ensuremath{\bold T^{#1}_{#2}}}
\providecommand{\lrb}[1]{\ensuremath{\left(#1\right)}}
\providecommand{\abs}[1]{\left\lvert #1\right\rvert}

\providecommand{\af}[1]{\text{A}^{\geq 0}(#1)}
\providecommand{\afs}[1]{\text{A}^{> 0}(#1)}
\DeclareMathOperator{\spec}{Spec}
\author{Brett Parker   }
\email{brettdparker@gmail.com}  
\thanks{This paper grew out of conversations I had with Bumsig Kim while visiting the Korean Institute for Advanced Studies. It was written during my stay at the Max Plank Institute for Mathematics in Bonn}
  
\title{Log geometry and exploded manifolds}

\begin{document}
\maketitle

\begin{abstract}Log Gromov-Witten invariants have recently been defined separately by Gross and Siebert and Abramovich and Chen. This paper provides a dictionary between log geometry and holomorphic exploded manifolds in order to compare Gromov-Witten invariants defined using exploded manifolds or log schemes. The gluing formula for Gromov-Witten invariants of exploded manifolds suggests an approach to proving analogous gluing formulas for log Gromov-Witten invariants.
\end{abstract}
\section{Introduction}
Around the turn of the century, Siebert suggested defining log Gromov-Witten invariants. Siebert's program has recently been carried out by Gross and Siebert \cite{GSlogGW}, and separately by Abramovich and Chen in \cite{acgw}, \cite{Chen}, \cite{Chen2} and \cite{acev}. Each of these groups built on work of Olsson \cite{olsson} which defined an appropriate deformation theory of log schemes, which has been used by Kim in \cite{kim} to define an obstruction theory for a stack of log curves. 

Meanwhile, over in the symplectic world, I have been working on an analogous project using exploded manifolds instead of log schemes. Gromov-Witten invariants of exploded manifolds are defined and  the associated gluing theorems are proved in \cite{egw}. I believe that each of these programs will turn out to define identical invariants in the settings to which they mutually apply (such as defining Gromov-Witten invariants relative to simple normal crossing divisors.) 

\

In section \ref{el}, some holomorphic exploded manifolds are described as log schemes. Then, in section \ref{expl} some log schemes will return the favor, and be describable in terms of exploded manifolds using a functor called the explosion functor. These explodable log schemes that correspond to exploded manifolds include log smooth schemes with characteristic or ghost sheaf $\mathcal M_{X}/\mathcal O^{*}_{X}$ a sheaf of toric monoids. In particular, the log smooth Deligne-Faltings log schemes used in \cite{acgw}  are explodable. 

 The relationship between the moduli stack of exploded curves and the moduli stack of basic log curves  is described in section \ref{moduli stack}. In particular, to obtain the moduli stack of curves in the explosion of $X^{\dag}$, the explosion functor is applied to the moduli stack of curves in $X^{\dag}$.

In section \ref{curves} the explosion functor is applied to log curves to obtain families of curves in the category of exploded manifolds. Section \ref{tropical} contains a comparison of the tropicalization defined by Gross and Siebert in \cite{GSlogGW}  and the tropical part functor defined for exploded manifolds. The minimality condition for curves used by Abramovich and Chen in \cite{acgw} and the basic condition for curves  defined by Gross and Siebert in \cite{GSlogGW} is translated to a condition for families of exploded curves in section \ref{basic}. Section \ref{cotangent} contains the observation that the cotangent sheaf of an exploded manifold is just the relative cotangent sheaf of a log scheme, then in section \ref{families}, the notion of a family of exploded manifolds is translated into the language of log geometry. Section \ref{reverse expl} is a list of log equivalents for exploded manifolds and maps - this includes subsection \ref{glog} which defines a generalized kind of log scheme to correspond to the exploded manifolds which are not representable as log schemes.

\

Sections \ref{refinements} and \ref{gluing} contain suggestions for the development of log Gromov-Witten invariants. Section \ref{refinements}  translates the notion of a refinement into the language of log schemes. This is important because given the correct cohomology theory, Gromov-Witten invariants are invariant under the operation of refinement. The moduli stack of curves in a refinement of $X^{\dag}$ should be a refinement of the moduli stack of curves in $X^{\dag}$. 

In section \ref{gluing}, I suggest an approach to a gluing formula for log Gromov-Witten invariants which is intuitive from the perspective of exploded manifolds, but may not be obvious from the perspective of log geometry.

\

In section \ref{smooth log}, I suggest a $C^{\infty}$ analogue of log schemes in which the construction of log Gromov-Witten invariants should be able to be imitated. Section \ref{dbarlog ncd} contains a definition of normal crossing divisors in this setting.

\section{Conventions for log schemes}

At this point I must admit an embarrassing lack of knowledge of algebraic geometry, so I shall only work with schemes $X$ over $\mathbb C$, and use the analytic topology. A point in $X$ will always mean a map $\spec \mathbb C\longrightarrow X$.  

I intend this paper to be read by an audience more familiar with log geometry than I am, but will include some definitions for the curious reader who does not know any log geometry. (The standard reference is \cite{kk}.) 
A log structure on a scheme $X$ is a sheaf $\mathcal M_{X}$ of commutative monoids on $X$ with a map
\[\alpha:\mathcal M_{X}\longrightarrow \mathcal O_{X} \] 
which is an isomorphism restricted to the inverse image of $\mathcal O^{*}_{X}\subset\mathcal O_{X}$. This map should be considered as a map of sheaves of monoids on $X$, where the monoid operation on  $\mathcal O_{X}$ is multiplication of functions. Denote $X$ with its log structure by $(X,\mathcal M_{X})$. Use the notation $\bar{\mathcal M}_{X}$ to denote the quotient $\mathcal M_{X}/\mathcal O^{*}_{X}$. I shall not use  the notation $\totb{X}$ for the underlying scheme of a log scheme, because I use that notation for the tropical part of an exploded manifold.

\section{Holomorphic exploded manifolds as log schemes}\label{el}

In this section, some holomorphic exploded manifolds are described as log schemes. Further intuition for understanding exploded manifolds in terms of log schemes will come from the explosion functor described in section \ref{expl}.

\

In \cite{iec}, it is proved that a full subcategory of holomorphic exploded manifolds may be regarded as log schemes over a point \[p^{\dag}=(\spec \mathbb C,\mathbb C^{*}\e{[0,\infty)})\] with a special log structure. 
As a monoid, $\mathbb C^{*}\e{[0,\infty)}$ is the direct sum of $\mathbb C^{*}$ with the operation of multiplication and $[0,\infty)$ with the operation of addition. The  structure map 
\[\alpha:\mathbb C^{*}\e{[0,\infty)}\longrightarrow \mathbb C\]
is referred to as the smooth part homomorphism in \cite{iec}.
\begin{equation}\label{smooth part def}\alpha(c\e x):=\totl{c\e x}:=\begin{cases}c\text{ if }x=0
\\ 0\text { if }x>0\end{cases}\end{equation}

\

\

A holomorphic exploded manifold $\ex B$ is defined as a topological space $\ex B$ with a sheaf $\mathcal E^{\times}(\ex B)$ of $\mathbb C^{*}\e {\mathbb R}$ valued functions locally isomorphic to an open subset of $\et mP$.
Accordingly, if $\ex B$ may be described as a log scheme, it is a log scheme \[\ex B^{\dag}:= (\totl{\ex B},\mathcal M_{\ex B})\longrightarrow p^{\dag}\] locally isomorphic to an open subset of some model log scheme \[(\et mP)^{\dag}:=(\totl{\et mP},\mathcal M_{\et mP})\longrightarrow p^{\dag}\]

To describe exploded manifolds in terms of log schemes, we must describe these models $(\et mP)^{\dag}$. As the information in the exploded manifold $\et mP$ and the log scheme $(\et mP)^{\dag}$ is identical, we shall use the notation $\et mP$ for both where no confusion shall arise.

\

The notation  $P$ above stands for an $m$ dimensional integral affine polytope in $\mathbb R^{m}$ - in other words $P$ is an $m$ dimensional subset of $\mathbb R^{m}$ which is the solution to a finite number of integral affine inequalities (which may be strict or non strict).  $\et mP$ may be described as a log scheme if and only if $P$ contains no entire lines. 
\begin{defn}
An integral affine map $\mathbb R^{m}\longrightarrow \mathbb R^{n}$ is a map in the form of \[x\mapsto Ax+y\] where $A$ is an integer matrix and $y\in\mathbb R^{n}$. 

\

Let $\af P$ indicate the monoid of integral affine maps \[P\longrightarrow [0,\infty)\]
 and $\afs P$ indicate the monoid of integral affine maps \[P\longrightarrow (0,\infty)\] 
\end{defn}

\

\
\begin{defn}For an integral affine polytope $P\subset\mathbb R^{m}$, define the scheme $\totl{\et mP}$ to be 
\[\totl{\et mP}:=\spec\lrb{\mathbb C[\af P]/\mathbb C[\afs P]}\]

where  $\mathbb C[\af P]$ indicates the monoid ring on $\af P$ over $\mathbb C$, and $\mathbb C[\afs P]$ indicates the ideal in this ring generated by $\afs P$. 
\end{defn}

Because $P$ is described by a finite number of inequalities,   $\mathbb C[\af P]/\mathbb C[\afs P]$ is always  finitely generated.

\

For example,  
\begin{itemize}
\item $\totl{\et m{[0,\infty)^{m}}}=\mathbb C^{m}$. 

\

\item If $P$ is a closed cone, then $\totl{\et mP}$ is the toric partial compactification of $(\mathbb C^{*})^{m}$ with fan $P$.

\

\item If $P$ is a standard closed $m$ dimensional simplex  $\totl{\et mP}$ is equal to 
\[\{z_{1}\dotsb z_{m+1}=0\}\subset\mathbb C^{m+1}\]

\

\item If $P$ is open, then $\totl{\et mP}$ is a point.
\end{itemize}

\

\begin{defn} The log scheme $(\et mP)^{\dag}$ is the scheme $\totl{\et mP}$ with the canonical log structure  associated to the map

\[\iota:\af P\longrightarrow \mathbb C[\af P]/\mathbb C[\afs P]\]

The definition of this canonical log structure is as in example 1.5(3) in \cite{kk}.

\

 In particular, $\af P$ considered as a constant sheaf on $\totl{\et mP}$ gives a pre log structure with structure map $\iota: \af P\longrightarrow \mathcal O_{\totl{\et mP}}$. The log structure on $(\et mP)^{\dag}$ is the log structure associated to this pre-log structure, so $\mathcal M_{\et mP}$ is constructed as a pushout in the category of sheaves of monoids on $\totl{\et mP}$ using the diagram
\[\begin{array}{ccc}\mathcal M_{\et mP}&\longleftarrow &\mathcal O^{*}_{\totl{\et mP}}
\\ \uparrow &&\uparrow
\\ \af P&\longleftarrow &\iota^{-1}\lrb{\mathcal O^{*}_{\totl{\et mP}}}
\end{array}\]

\end{defn}

 Roughly speaking, the sheaf of monoids $\mathcal M_{\et mP}$ on $\totl{\et mP}$ is obtained by taking the direct sum  of $\mathcal O^{*}_{\totl{\et mP}}$ with the constant sheaf $\af P$, then quotienting by the equivalence generated by equating $(f,x)$ with $ (f\iota (x),0)$ if $\iota (x)$ is invertible. The map $\alpha:\mathcal M_{\et mP}\longrightarrow \totl{\et mP}$ is given by $(f,x)\mapsto f\iota(x)$.

\

\begin{remark}For the definitions within this paper, $m$ is the dimension of $P$ so the $m$ in $\et mP$ is redundant. Although it is not done in this paper, it would also be consistent to allow $P$ to be a polytope in $\mathbb R^{m}$ which has dimension less than $m$, and use $\af P$ to mean the integral affine functions on $\mathbb R^{m}$ which are non negative when restricted to $ P$. Then $\et mP=(\mathbb C^{*})^{m-n}\times\et nP$ when $P$ has dimension $n$.
\end{remark}

Note the following:\begin{itemize}
\item A point in $\totl{\et mP}$ corresponds to a homomorphism \[h:(\af P,+)\longrightarrow (\mathbb C,\times)\] which sends $\afs P$ to $0$.
\begin{itemize}\item The stalk of $\bar{\mathcal M}_{\et mP}$ at this point is $\af P/h^{-1}(\mathbb C^{*})$. \item The submonoid $h^{-1}(\mathbb C^{*})\subset \af P$ is some set of integral affine functions on $P$ which must have a common zero because if they did not, the sum of two of these functions would be in $\afs P$, and therefore not in $h^{-1}(\mathbb C^{*})\subset \af P$. \item Call the set of common zeros of a submonoid of $\af P$ a face of $P$. In fact, $h^{-1}(\mathbb C^{*})$ is equal to the submonoid of integral affine functions which vanish on some face $F$ of $P$, so the stalk of $\bar{\mathcal M}_{\et mP}$ at this point is equal to $\af F$.  
\end{itemize}

\

\item 
There is a canonical  injective map of the constant sheaf $\mathbb C^{*}\e{[0,\infty)}$ into $\mathcal M_{\et mP}$ by $c\e x\mapsto (c,x)$, where $c$ is regarded as a constant $\mathbb C^{*}$ valued function on $\totl{\et mP}$ and $x$ is a constant $\mathbb R$ valued function on $P$. This makes $\et mP$ into a log scheme over $p^{\dag}$.  

\

\item $\af P$ is a saturated sub monoid of the group of integral affine maps $P\longrightarrow \mathbb R$, so $\af P$ is a saturated integral monoid. Similarly, the stalk of  $\bar{\mathcal M}_{\et mP}$ at any point is $\af F$ for some face $F$ of $P$ so $\et mP$ is a saturated integral log scheme.

\
 
\item The description of $\ex B$ as a log scheme locally isomorphic to some subset of $\et mP$ implies that $\ex B$ is a quasi coherent log scheme as defined in \cite{kk}. As $[0,\infty)$ is not a finitely generated monoid, exploded manifolds are not coherent log schemes.  On the other hand, $\af P$ is finitely generated as a monoid over the constant maps to $[0,\infty)$, so exploded manifolds are as close to being coherent log schemes as possible.

\

\item If $P$ is a closed cone, then $\af P$ decomposes as a direct sum of $[0,\infty)$ with the monoid $Q$ of linear integral maps $P\longrightarrow[0,\infty)$. In this case, 
\[\mathbb C[Q]=\mathbb C[\af P]/\mathbb C[\afs P]\]
The spectrum of $\mathbb C[Q]$ is the toric space with fan $P$, and the log structure associated to 
\[Q\longrightarrow \mathbb C[Q]\]
is the standard log structure on this toric space with sheaf of monoids given by the functions which are invertible when restricted to $(\mathbb C^{*})^{m}$ inside of the spectrum of $\mathbb C[Q]$. The sheaf $\mathcal M_{\et mP}$ is simply the direct sum of this standard log structure with the constant sheaf $[0,\infty)$.   This example shall be revisited when we consider the explosion functor.
\end{itemize}

\

\

An exploded manifold $\ex B$ may be recovered from $\ex B^{\dag}\longrightarrow p^{\dag}$ as follows:
\begin{itemize}
\item The set of points in the exploded manifold $\ex B$ is equal to the set of  maps   $p^{\dag}\longrightarrow \ex B^{\dag}$. In other words, if we consider $p^{\dag}$ as a log scheme over itself, $\ex B$ as a set is the set of maps of $p^{\dag}$ into $\ex B^{\dag}$ as a log scheme over $p^{\dag}$.

 For example, a point in $\et mP$ is equivalent to a homomorphism \[\af P\longrightarrow \mathbb C^*\e{[0,\infty)}\] sending the constant function $y$ to  $1\e y$,  so that the smooth part of the image of any function in $\afs P$ is $0$. Such a homomorphism may always be written uniquely in the form \[(n\cdot x+y)\mapsto \e y\prod (c_{i}\e {a_{i}})^{n_{i}}\] where $(c_{1}\e {a_{1}},\dotsc,c_{m}\e {a_{m}} )\in (\mathbb C^{*}\e{\mathbb R})^{m}$ is so that $(a_{1},\dotsc, a_{m})\in P$.
\item A map $p^{\dag}\longrightarrow \ex B$ includes the information of a map $\spec \mathbb C\longrightarrow \totl{\ex B}$, so there is a set map  \[\ex B\longrightarrow \totl{\ex B}\]
The topology on $\ex B$ is the topology pulled back from the analytic topology on $\totl{\ex B}$.
\item It is easy to check that this map $\ex B\longrightarrow \totl{\ex B}$ is surjective, so we may regard $\mathcal M_{\ex B}$ as a sheaf on $\ex B$. 
\item In terms of exploded manifolds, $\mathcal M_{\ex B}$ actually the sheaf of holomorphic maps to $\et 1{[0,\infty)}$. We may identify the set of points in $\et 1{[0,\infty)}$ with $\mathbb C^{*}\e{[0,\infty)}$  so that $\mathcal M_{\ex B}$ may be regarded as a sheaf of $\mathbb C^{*}\e{[0,\infty)}$ valued functions. The map $\alpha: \mathcal M_{\ex B}\longrightarrow\mathcal O_{\totl{\ex B}}$ is given by composition with  the smooth part homomorphism defined on page \pageref{smooth part def} \[\alpha(f):=\totl{f}\]
\item The defining sheaf $\mathcal E^{\times }(\ex B)$ of the exploded manifold $\ex B$ is  equal to the sheaf $(\mathcal M_{\ex B})^{gr}$ of groups generated by $\mathcal M_{\ex B}$, which can also be regarded as a sheaf of $\mathbb C^{*}\e{\mathbb R}$ valued functions. 
\item There is an exploded manifold $\ex T:=\et 1{\mathbb R}$ so that $\mathcal E^{\times}(\ex B)$ is the sheaf of holomorphic maps to $\ex T$. This exploded manifold $\ex T$ may not be described as a log scheme. 

\end{itemize}

 The holomorphic exploded manifolds $\ex B$ which may be regarded as log schemes $\ex B^{\dag}$ are those for which $\mathcal E^{\times }(\ex B)$ is generated by functions taking values in $\mathbb C^{*}\e{[0,\infty)}\subset\mathbb C^{*}\e{\mathbb R}$ - in this case $\mathcal M_{\ex B}$ is simply the subsheaf of $\mathcal E^{\times}(\ex B)$ consisting of those functions taking values in $\mathbb C^{*}\e{[0,\infty)}\subset\mathbb C^{*}\e{\mathbb R}$.

\

For general integral affine polytopes, $\et mP$ is not a log scheme. We may however describe all maps from $\ex B$ to $\et mP$ if $\ex B$ is describable as a log scheme $\ex B^{\dag}$. A map $\ex B^{\dag}\longrightarrow \et mP$ is equivalent to $m$ global sections of $(\mathcal M_{\ex B})^{gr}$, which evaluated at any point in $\ex B^{\dag}$ give $(c_{1}\e {a_{1}},\dotsc,c_{m}\e{a_{m}})$ so that  $(a_{1},\dotsc,a_{m})\in P$. We may rephrase this information as follows: 

A map $\ex B\longrightarrow \et mP$ is equivalent to a homomorphism from the group of integral affine functions on $P$ to the group of global sections of $(\mathcal M_{\ex B})^{gr}$ which sends $\af P$ to $\mathcal M_{\ex B}$, sends $\afs P$ to $\alpha^{-1}(0)\subset\mathcal M_{\ex B}$, and which is the canonical inclusion of $[0,\infty)$ in $\mathcal M_{\ex B}$ when restricted to constant functions in $\af P$.

\section{The explosion functor }\label{expl}

The explosion functor allows us to think of some log schemes as exploded manifolds. In particular, given a log scheme $(X,\mathcal M_{X})$ over $\mathbb C$, the explosion of $(X,\mathcal M_{X})$ is given by a base change to a log scheme over $p^{\dag}$

\[\begin{array}{ccc}\expl(X,\mathcal M_{X})&\longrightarrow &(X,\mathcal M_{X})
\\ \downarrow && \downarrow
\\ p^{\dag}&\longrightarrow &(\spec\mathbb C,\mathbb C^{*})\end{array}\]

\

\

Explicitly, 
\[\expl(X,\mathcal M_{X}):=(X,\mathcal M_{X}\oplus[0,\infty))\]
where the new structure map
\[\alpha':\mathcal M_{X}\oplus[0,\infty)\longrightarrow \mathcal O_{X}\] 
is described in terms of the old structure map
\[\alpha:\mathcal M_{X}\longrightarrow \mathcal O_{X}\]
by 
\[\alpha'(f,x):=\begin{cases}\alpha(f)\text{ if }x=0
\\ 0\text{ if }x>0
\end{cases}\]
The map $\expl(X,\mathcal M_{X})\longrightarrow p^{\dag}$ is given by extending the inclusion of the constant sheaf $\mathbb C^{*}\hookrightarrow \mathcal M_{X}$ to the obvious inclusion $\mathbb C^{*}\oplus[0,\infty)\hookrightarrow \mathcal M_{X}\oplus [0,\infty)$. 

\

For example, 

\begin{enumerate}
\item Let $Q$ be a toric monoid consisting of the nonnegative integral linear functions on some closed $m$ dimensional integral affine cone $P$ which contains no lines. Consider the  log structure on $\spec \mathbb C[Q]$ given by the map 
\[Q\longrightarrow \mathbb C[Q]\]
The explosion of this log scheme is $\et mP$. 

\

\item Let $Q$ be the above toric monoid. There is a log point $(\spec \mathbb C,\mathbb C^{*}\oplus Q)$ with structure map
\[\alpha(c,q):=\begin{cases}c\text{ if }q=0
\\ 0\text{ if }q\neq 0\end{cases}\]
The explosion of this log point is 
\[\expl(\spec \mathbb C,\mathbb C^{*}\oplus Q)=\et m{P^{\circ}}\] where $P^{\circ}$ is the interior of $P$.

\

\item Let $Q$ be the above toric monoid, and let $P'$ be an integral affine cone with closure equal to $P$. Then let $Q'\subset Q$ be the submonoid consisting of integral linear functions which are strictly positive on $P'$. Let $X_{P'}^{\dag}$  be the log scheme  given by the map 
\[Q\longrightarrow \mathbb C[Q]/\mathbb C[Q']\]
Then
\[\expl X_{P'}^{\dag}=\et m{P'}\]

\

\item
Let \[f:P_{1}\longrightarrow P_{2}\] be a surjective map of integral affine cones so that every face of $P_{1}$ maps surjectively onto a face of $P_{2}$ and so that the corresponding homomorphism
\[f^{*}:Q_{2}\longrightarrow Q_{1}\] 
is injective and the corresponding map of groups generated by $Q_{i}$ has torsion free cokernel.
There is a corresponding map of  the log spaces
\[f^{\dag}:X^{\dag}_{P_{1}}\longrightarrow X^{\dag}_{P_{2}}\]
defined as in the above example.
The explosion of $f^{\dag}$
\[\expl f^{\dag}:\et {m_{1}}{P_{1}}\longrightarrow \et {m_{2}}{P_{2}}\]
is a submersion. A point $p^{\dag}\longrightarrow \et {m_{2}}{P_{2}}$ is equivalent to the choice of a homomorphism $Q_{2}\longrightarrow \mathbb C^{*}\e{[0,\infty)}$ which sends the strictly positive integral linear functions to $\mathbb C^{*}\e{(0,\infty)}$. In particular, this includes the choice of a homomorphism $Q_{2}\longrightarrow [0,\infty)$ which is strictly positive on the strictly positive functions in $Q_{2}$. This is equivalent to the choice of a point in $P_{2}$. Suppose that this point is in the interior of $P_{2}$, and let $P'$ indicate inverse image of this point under $f$. Then the fiber product of $\expl f^{\dag}$ with $p^{\dag}\longrightarrow \et {m_{2}}{P_{2}}$ is isomorphic to $\et {m_{1}-m_{2}}{P'}$. 
\[\begin{array}{ccl}\et {m_{1}-m_{2}}{P'}&\longrightarrow &\et {m_{1}}{P_{1}}
\\ \downarrow&&\downarrow \expl f^{\dag}
\\   p^{\dag}&\longrightarrow&\et {m_{2}}{P_{2}} \end{array}\]

\end{enumerate}

\

\

The explosion of any log scheme locally isomorphic to an open subset of one of first three of the above examples is an exploded manifold. 

\begin{defn}\label{xdef} A log scheme $(X,\mathcal M_{X})$ over $\mathbb C$ is explodable if it is locally isomorphic to  an open subset of $X^{\dag}_{P}$ where $P\subset[0,\infty)^{m}$ is an integral affine cone.

\

 $X^{\dag}_{P}$ is defined as follows: Let $Q$ be the toric monoid consisting of integral linear maps $P\longrightarrow [0,\infty)$ and $Q'\subset Q$ be the submonoid consisting of integral linear maps $P\longrightarrow (0,\infty)$. Then $X^{\dag}_{P}$ is the log scheme associated to the map 
 \[Q\longrightarrow \mathbb C[Q]/\mathbb C[Q']\]
where $\mathbb C[Q']$ indicates the ideal in $\mathbb C[Q]$ which is generated by $Q'$ if $Q'$ is non-empty, and the zero ideal if $Q'$ is empty.
\end{defn}

For example, the log structure associated to a normal crossing divisor is explodable, and any log smooth scheme $(X,\mathcal M_{X})$ with  $\bar{\mathcal M}_{X}$ a sheaf of toric monoids is also explodable. In particular, the log smooth Deligne-Faltings log structures used by Abramovich and Chen in \cite{acgw} and Chen in \cite{Chen} and \cite{Chen2} are explodable. The explosion functor applied to an explodable log scheme is a holomorphic exploded manifold.

\

Note that  $P$ being a closed integral affine cone  is equivalent to $X_{P}^{\dag}$ being log smooth.

\

The following are some further examples to give intuition about what the explosion functor does. 

\begin{enumerate}\item A complex algebraic manifold $M$ can be regarded as a log scheme with the trivial log structure, $(M,\mathcal O^{*}_{M})$. Applying the explosion functor gives the log scheme $(M,\mathcal O^{*}_{M}\oplus[0,\infty))\longrightarrow (\spec \mathbb C,\mathbb C^{*}\oplus[0,\infty))$. Let $M_{1}$ and $M_{2}$ be complex manifolds. A map 
\[\expl M_{1}\longrightarrow \expl M_{2}\]
is equivalent to a map 
\[f:M_{1}\longrightarrow M_{2}\] 
and a homomorphism 
\[h:f^{-1}\mathcal O^{*}_{M_{2}}\oplus[0,\infty)\longrightarrow \mathcal O^{*}_{M_{1}}\oplus[0,\infty)\]
compatible with $f$ and the structure maps, and also compatible with the maps to $p^{\dag}$. Compatibility with the inclusions  of $\mathbb C^{*}\oplus[0,\infty)$ coming from the maps to $p^{\dag}$ gives that $h$ must be the identity homomorphism on $[0,\infty)$. The following diagram must commute, 
\[\begin{array}{ccc}f^{-1}\mathcal O^{*}_{ M_{2}}\oplus[0,\infty)&\longrightarrow &\mathcal O^{*}_{M_{1}}\oplus[0,\infty)\\\downarrow &&\downarrow
\\ f^{-1}\mathcal O_{M_{2}}&\xrightarrow {\circ f}&\mathcal O_{M_{1}}
\end{array}\]
so the unique homomorphism $h$ satisfying the required conditions is
\[h(g,x):=(g\circ f,x)\]

This implies that the explosion functor identifies the category of  complex manifolds as a full subcategory of the category of holomorphic  exploded manifolds.

\

\item A point in $\et mP$ corresponds to a map $p^{\dag}\longrightarrow \expl X^{\dag}_{P}$. Because $\expl X^{\dag}_{P}$ is given as a base change,

\[\begin{array}{ccc}\expl X^{\dag}_{P}&\longrightarrow& X^{\dag}_{P}
\\ \downarrow&&\downarrow
\\ p^{\dag}&\longrightarrow& (\spec \mathbb C,\mathbb C^{*})\end{array}\]
a map from $p^{\dag} $ to $\expl X^{\dag}_{P}$ is equivalent to a map from $p^{\dag}$ to $X^{\dag}_{P}$. 

Recall that $X^{\dag}_{P}$ is defined as the log scheme associated to the homomorphism
 \[Q\longrightarrow \mathbb C[Q]/\mathbb C[Q']\]
where $Q$ is the toric monoid given by nonnegative integral linear functions on $P$ and $Q'$ is the submonoid consisting of strictly positive integral linear functions. It follows that a map $p^{\dag}\longrightarrow X^{\dag}_{P}$ is equivalent to a commutative diagram of  homomorphisms 
\[\begin{array}{ccl}Q&\longrightarrow &\mathbb C^{*}\e{[0,\infty)}
\\ \downarrow &&\downarrow\totl{\cdot}
\\\mathbb C[Q]/\mathbb C[Q']&\longrightarrow &\mathbb C\end{array}\]
The homomorphism $Q\longrightarrow \mathbb C^{*}\e{[0,\infty)}$ gives in particular a homomorphism from $Q$ to $[0,\infty)$. As $Q$ is given by  the nonnegative integral linear functions on $P$, this homomorphism corresponds to a point in the closure of $P$. The requirement that the homomorphism $Q\longrightarrow \mathbb C^{*}\e{[0,\infty)}$ composed with the smooth part homomorphism vanishes on $Q'$ is equivalent to requiring that that this point must be inside $P$. The homomorphism $Q\longrightarrow \mathbb C^{*}$ may be chosen independent of this homomorphism to $[0,\infty)$. The group generated by $Q$ is $\mathbb Z^{m}$, so a homomorphism $Q\longrightarrow \mathbb C^{*}$ is equivalent to a point in $(\mathbb C^{*})^{m}$. A point  in $T^{m}_{P}=\expl X^{\dag}_{P}$ is therefore equivalent to a point in $P$ and a point in $(\mathbb C^{*})^{m}$. More generally, a map from (the explosion of) a complex manifold to $\expl X^{\dag}_{P}$ is equivalent to a map to $(\mathbb C^{*})^{m}\times P$ where $P$ is given the discrete topology. 
The points of $\expl X^{\dag}_{P}$ that correspond to points in $X^{\dag}_{P}$ (ie maps $\spec \mathbb C\longrightarrow X^{\dag}_{P}$) are the points over $0\in P$.

\

\item In the notation of definition \ref{xdef}, $X^{\dag}_{[0,\infty)}$ is the log space associated to the homomorphism
\[\mathbb N\longrightarrow\mathbb C[\mathbb N]\]
which is also the log scheme associated to the divisor $0$ in $\mathbb C$. 

A map from any log scheme $(X,\mathcal M_{X})$ to $X^{\dag}_{[0,\infty)}$ is equivalent to a global section of $\mathcal M_{X}$.
Accordingly, a map from a  complex manifold $M$ to $X^{\dag}_{[0,\infty)}$ is equivalent to a map  $M\longrightarrow \mathbb C^{*}\subset X^{\dag}_{[0,\infty)}$. 

As $\expl X^{\dag}_{[0,\infty)}$ is just given by a base change, a map $\expl M\longrightarrow \expl X^{\dag}_{[0,\infty)}$ is  equivalent  to a pair of maps 
\[\begin{array}{ccc}\expl M&\longrightarrow &X^{\dag}_{[0,\infty)}
\\\downarrow 
\\ p^{\dag}\end{array}\]
so a map $\expl M\longrightarrow \expl X^{\dag}_{[0,\infty)}$  is equivalent to a choice of global section of $\mathcal M_{\expl M}:=\mathcal O^{*}_{M}\oplus [0,\infty)$. If $M$ is connected, this corresponds to a nonvanishing holomorphic function and a constant $x\in[0,\infty)$. If $x=0$, these maps correspond to the maps $M\longrightarrow X^{\dag}_{[0,\infty)}$, but if $x>0$, then these maps do not correspond to any map $M\longrightarrow X^{\dag}_{[0,\infty)}$ because their image in the underlying scheme of $X^{\dag}_{[0,\infty)}$ is $0\in\mathbb C$, where there is a nontrivial log structure.

\

 Even though there is no map  $M\longrightarrow X^{\dag}_{[0,\infty)}$ with image  the special point $0$ in the underlying scheme $\mathbb C$ of $X^{\dag}_{[0,\infty)}$, there are maps which could be considered as families of maps from $M$ to $X^{\dag}_{[0,\infty)}$ with image $0\in\mathbb C$. Make a base change: 
\[\begin{array}{ccc}M^{\dag}&\longrightarrow& M
\\\downarrow && \downarrow
\\X^{\dag}_{(0,\infty)}:=(\spec \mathbb C,\mathbb C^{*}\oplus \mathbb N)&\longrightarrow &(\spec \mathbb C,\mathbb C^{*})\end{array}\]
so $M^{\dag}=(M,\mathcal O^{*}_{M}\oplus \mathbb N)$. Again, a map $M^{\dag}\longrightarrow X^{\dag}_{[0,\infty)}$ is equivalent to a global section of $\mathcal O^{*}_{M}\oplus\mathbb  N$, which is equivalent to a nonvanishing function and a natural number. If this number is positive,  the image of this map is $0\in\mathbb C$. 

A map $M^{\dag}\longrightarrow X^{\dag}_{[0,\infty)}$ may be thought of in terms of maps of exploded manifolds as follows:
The map $M^{\dag}\longrightarrow X^{\dag}_{(0,\infty)}$ is like a trivial family with fiber equal to $M$. Its explosion $\expl M^{\dag}\longrightarrow \expl X^{\dag}_{(0,\infty)}$ actually is a trivial family of exploded manifolds with each fiber equal to $\expl M$. Correspondingly, 
the explosion of a map $M^{\dag}\longrightarrow X^{\dag}_{[0,\infty)}$ is a family of maps from $\expl M$ to $\expl X^{\dag}$.
\[\begin{array}{cc}\expl M^{\dag}&\longrightarrow \expl X^{\dag}_{[0,\infty)}
\\\downarrow 
\\  \expl X^{\dag}_{(0,\infty)}&\end{array}\]
A point in $X^{\dag}_{(0,\infty)}$ corresponds to an element $c\e x\in\mathbb C\e{(0,\infty)}$. Restricting the map corresponding to the section $(f,n)$ of $\mathcal M_{M^{\dag}}$ to this point gives the map $\expl M\longrightarrow \expl X^{\dag}_{[0,\infty)}$ corresponding to the section $(c^{n}f,nx)$.

\end{enumerate}

\section{curves}\label{curves}

In this section, we compare curves in the category of exploded manifolds to the explosion functor applied to the log smooth curves found in \cite{GSlogGW} or \cite{Chen}. In particular, consider a prestable log curve 

\[\begin{array}{c}(C,\mathcal M_{C})
\\\downarrow
\\ (\spec\mathbb C,\mathbb C^{*}\oplus Q)\end{array}\]
where $Q$ is a toric monoid. In the notation of definition \ref{xdef}, \[(\spec\mathbb C,\mathbb C^{*}\oplus Q)=X^{\dag}_{P}\] where $P$ is the open cone which is the space of homomorphisms from $Q$ to $\mathbb R$ which are strictly positive on nonzero elements. If $Q\neq 0$, then the explosion of this curve will always be considered as a family of exploded curves rather than an individual exploded curve.

\

Locally, $(C,\mathcal M_{C})\longrightarrow X^{\dag}_{P}$ is isomorphic to an open subset of one of the following $3$  models:
\begin{enumerate}
\item At nodes of $C$, $(C,\mathcal M_{C})\longrightarrow X^{\dag}_{P}$ is isomorphic to an open subset of the log scheme induced from a homomorphism
\[Q\oplus_{\mathbb N}\mathbb N^{2}\longrightarrow\mathbb C[z,w]/zw\]
\[(q,a,b)\mapsto\begin{cases}z^{a}w^{b}\text{ if }q=0
\\ 0 \text{ if }q\neq 0\end{cases}\]
where $Q\oplus_{\mathbb N}\mathbb N^{2}$ is defined by the pushout diagram
\[\begin{array}{ccc}Q\oplus_{\mathbb N}\mathbb N^{2}&\longleftarrow &\mathbb N^{2}
\\ \uparrow &&\uparrow
\\ Q&\longleftarrow &\mathbb N
\end{array}\]
with the map $\mathbb N\longrightarrow\mathbb N^{2}$ the diagonal map, and the homomorphism $\mathbb N\longrightarrow Q$ nonzero, and determined by a nonzero element $\rho_{q}\in Q$ which is the image of $1$.
The map from this log scheme to $X^{\dag}_{P}$ is determined by the homomorphism $Q\longrightarrow Q\oplus_{\mathbb N}\mathbb N^{2}$.

\

In terms of the notation from definition \ref{xdef}, $(C,\mathcal M_{C})\longrightarrow X^{\dag}_{P}$ is locally in the form of an open subset of a map $X^{\dag}_{P'}\longrightarrow X^{\dag}_{P}$ corresponding to some integral linear map $P'\longrightarrow P$ given as follows:

If $P\subset \mathbb R^{m}$, $P'$ is the intersection of the inverse image of $P$ under the standard projection $\mathbb R^{m+1}\longrightarrow \mathbb R^{m}$ with two half spaces
\[P':=\{(x_{1},\dotsc,x_{m+1})\text{ so that } (x_{1},\dotsc,x_{m})\in P,\  x_{m+1}\geq 0,\ x_{m+1}\leq \rho_{q}\}\]
 where $\rho_{q}\in Q$ is identified with an integral linear function on $\mathbb R^{m}$ which pulls back to an integral linear function on $\mathbb R^{m+1}$. 
 
 Note that in this setting,  $Q\oplus_{\mathbb N}\mathbb N^{2}$ is the monoid of integral linear maps $P'\longrightarrow[0,\infty)$, and $ \mathbb C[x,y]/xy$ is equal to $\mathbb C[\af {P'}]/\mathbb C[\afs {P'}]$. The description of $Q\oplus_{\mathbb N}\mathbb N^{2}$ as a pushout corresponds to a description of $P'$  as a fiber product via the diagram
 \[\begin{array}{ccl}P'&\longrightarrow &[0,\infty)^{2}
 \\ \downarrow &&\ \ \ \downarrow x_{1}+x_{2}
 \\ P&\xrightarrow{\rho_{q}}&[0,\infty) \end{array}\]
 which corresponds to a description of $X^{\dag}_{P'}$ as a fiber product via the diagram
 \[\begin{array}{ccc}X^{\dag}_{P'}&\longrightarrow &X^{\dag}_{[0,\infty)^{2}}
 \\\downarrow &&\downarrow
 \\ X^{\dag}_{P}&\longrightarrow &X^{\dag}_{[0,\infty)}
\end{array} \]

\item At marked points, $(C,\mathcal M_{C})\longrightarrow X^{\dag}_{P}$ is locally in the form of an open subset of the log scheme corresponding to the homomorphism
\[Q\oplus \mathbb N\longrightarrow \spec \mathbb C[z]\]
\[(q,a)\mapsto \begin{cases}z^{a}\text{ if }q=0
\\ 0 \text{ if }q\neq 0\end{cases}\]
(The marked point in this local model is at $z=0$ where the log structure is special.) The map to $X^{\dag}_{P}$ is given by the obvious homomorphism $Q\longrightarrow Q\oplus\mathbb N$. 

In terms of the notation of definition \ref{xdef}, $(C,\mathcal M_{C})$ around a marked point is isomorphic to an open subset of a map $X^{\dag}_{P'}\longrightarrow X^{\dag}_{P}$ where $P'$ is the inverse image of $P$ under the projection $\mathbb R^{m+1}\longrightarrow \mathbb R^{m}$ intersected with the half space
\[x_{m+1}\geq 1\]
We may also describe $X^{\dag}_{P'}$ as a  product via the diagram
\[\begin{array}{ccc}X^{\dag}_{P'}&\longrightarrow&X^{\dag}_{[0,\infty)}
\\ \downarrow &&\downarrow
\\X^{\dag}_{P}&\longrightarrow&(\spec\mathbb C,\mathbb C^{*})\end{array}\]

\item Elsewhere, $(C,\mathcal M_{C})$ is locally in the form of an open subset of 
\[(\mathbb C,\mathcal O^{*}(\mathbb C)\oplus Q)\]
with the structure map 
\[\alpha(f,q):=\begin{cases}f\text{ if }q=0
\\ 0\text{ if }q\neq 0\end{cases}\]
and the obvious projection to $X^{\dag}_{P}$ coming from the identity map $Q\longrightarrow Q$. This is just the product of $X^{\dag}_{P}$ with $\mathbb C$. Sufficiently small open subsets of this model are also isomorphic to open subsets of the other two models where $\bar{\mathcal M}_{X^{\dag}_{P'}}$ is  the constant sheaf $Q$. 
\end{enumerate}

\

\

The prestable log curves defined by Kim in  \cite{kim} are similar, except they are locally isomorphic to the first and third models. 

\

As all three local models for a curve over $X^{\dag}_{P}$ are explodable, exploding $(C,\mathcal M_{C})\longrightarrow X^{\dag}_{P}$ gives a map of exploded manifolds, which turns out to be a family of exploded curves.
\[\begin{array}{c}\expl(C,\mathcal M_{C})
\\ \downarrow
\\ \expl X^{\dag}_{P}\end{array}\]

\

Choose a point in $\expl X^{\dag}_{P}$ (in other words a map $p^{\dag}\longrightarrow \expl X^{\dag}_{P}$). The fiber of $\expl(C,\mathcal M_{C})\longrightarrow \expl X^{\dag}_{P}$ over this point is a curve $\ex C$ in the category of exploded manifolds. 
\[\begin{array}{ccc}\ex C&\longrightarrow &\expl(C,\mathcal M_{C})
\\ \downarrow&&\downarrow
\\ p^{\dag}&\longrightarrow&\expl X^{\dag}_{P}\end{array}\]
Apart from a unique curve $\et 1{\mathbb R}$ which can not be described as a log scheme, all abstract curves in the category of exploded manifolds may be constructed as the explosion of some prestable log curve as above.

\

The underlying scheme $\totl{\ex C}$ of $\ex C$ is still the original nodal curve $C$. The description of $(C,\mathcal M_{C})$ locally in terms of fiber products translates to a description of $\ex C$ in terms of fiber products. In particular, $\ex C$ is locally isomorphic to one of the following $3$ models.
\begin{enumerate}
\item At nodes, $\expl(C,\mathcal M_{C})\longrightarrow \expl X^{\dag}_{P}:=\et mP$ is locally isomorphic to an open subset of $\et{m+1}{P'}\longrightarrow \et mP$ constructed via the pullback diagram
\[\begin{array}{ccc}\et {m+1}{P'}&\longrightarrow &\et 2{[0,\infty)^{2}}
\\ \downarrow &&\downarrow \tilde z_{1}\tilde z_{2}
\\ \et mP&\xrightarrow{\tilde z^{\rho_{q}}} &\et 1{[0,\infty)}\end{array}\]
The labeled maps above are labeled by the global sections of the defining sheaf of monoids using the notation  for these as $\mathbb C^{*}\e{[0,\infty)}$ valued maps defined in terms of coordinate functions $\tilde z$. In particular, $\tilde z^{\rho_{q}}$ corresponds to the integral linear function $\rho_{q}$, and $\tilde z_{1}\tilde z_{2}$ corresponds to the integral linear map $[0,\infty)^{2}\longrightarrow [0,\infty)$ given by $x_{1}+x_{2}$.

The function $\tilde z^{\rho_{q}}$ takes a value $c\e l\in\mathbb C^{*}\e{(0,\infty)}$ at $p^{\dag}\in X^{\dag}_{P}$. The fiber over $p^{\dag}$ is equal to the fiber over $c\e l\in \et 1{[0,\infty)}$, which is defined by the equation  $\tilde z_{1}\tilde z_{2}=c\e l$. This fiber is isomorphic to $\et 1{[0,l]}$, so $\ex C$ at a node is isomorphic to an open subset of $\et 1{[0,l]}$ for some positive $l$. The following is a pullback diagram
\[\begin{array}{ccccc}\et 1{[0,l]}&\longrightarrow &\et {m+1}{P'}&\longrightarrow &\et 2{[0,\infty)^{2}}
\\\downarrow&&\downarrow&&\downarrow
\\ p^{\dag}&\longrightarrow&\et mP&\longrightarrow &\et 1{[0,\infty)}\end{array}\]
Note that not all fibers of $\et 2{[0,\infty)^{2}}\longrightarrow \et 1{[0,\infty)}$ are isomorphic, as the isomorphism type of $\et 1{[0,l]}$ depends on $l$, and a fiber over $c\e 0$ is isomorphic to the complex manifold $\mathbb C^{*}$. How these charts (considered as open subsets of $\et 1{[0,l]}$) are glued to other charts on the curve $\ex C$ also depends on the full value of $\tilde z^{\rho_{q}}$ in $\mathbb C^{*}\e{[0,\infty)}$, not just on $l$.

\

\item At marked points, $\expl (C,\mathcal M_{C})\longrightarrow \et mP$ is locally isomorphic to an open subset of a product constructed via the pullback diagram
\[\begin{array}{ccc}\et mP\times \et 1{[0,\infty)}&\longrightarrow&\et 1{[0,\infty)}
\\ \downarrow &&\downarrow
\\ \et mP &\longrightarrow &p^{\dag}
\end{array}\]
The marked point in $C$ corresponds to $0\in \mathbb C=\totl{\et1 {[0,\infty)}}$.
The fiber over any point in $ \et mP$ is equal to $\et 1{[0,\infty)}$, so $\ex C$ is locally isomorphic to an open subset of $\et 1{[0,\infty)}$ at a marked point.
\item Elsewhere, $\expl (C,\mathcal M_{C})\longrightarrow \et mP$ is locally isomorphic to an open subset of $\et mP\times \expl \mathbb C\longrightarrow \et mP$. So $\ex C$ is locally isomorphic to an open subset of (the explosion of) $\mathbb C$.   
\end{enumerate} 

\

To define curves in the category of exploded manifolds, we need the concept of completeness. 

\begin{defn}A complete exploded manifold $\ex B$ is one for which $\totl{\ex B}$ is compact, and which is locally isomorphic to an open subset of $\et mP$ where $P$ is a closed integral affine polytope.
\end{defn}
 If $X^{\dag}$ is an explodable log scheme, $\expl X^{\dag}$ is complete if and only if $X^{\dag}$ is compact and log smooth. There is a similar notion of a complete map defined in \cite{iec}, which is not the same as being proper and log smooth. Completeness  is a property of exploded maps which plays the role of properness. 

\begin{defn}
A curve in the category of exploded manifolds is a one dimensional holomorphic exploded manifold $\ex C$ which is complete. 
\end{defn} In other words,  $\totl{\ex C}$ is compact and $\ex C$ is either $\et 1{\mathbb R}$ or locally isomorphic to an open subset of   $\et 1{[0,l]}$, $\et 1{[0,\infty)}$ or $\mathbb C$.

\begin{defn}A curve in a holomorphic exploded manifold $\ex B$ is a map $f:\ex C\longrightarrow \ex B$ where $\ex C$ is a curve. 
\end{defn}
There is a natural definition of stability for curves in exploded manifolds. If $\ex B$ may be regarded as a log scheme over $p^{\dag}$, then this stability condition translates to $\ex C$ being representable as a log scheme, and the underlying map $\totl f:\totl{\ex C}\longrightarrow \totl{\ex B}$ of schemes being stable. Therefore, if 
\[\begin{array}{ccc}&C^{\dag}&\xrightarrow{f^{\dag}} B^{\dag}
\\ &\downarrow &
\\ &X^{\dag}_{P} \end{array}\] is a curve over a point with explodable log structure in an explodable log scheme $B^{\dag}$, then 
\[\begin{array}{ccc}&\expl C^{\dag}&\xrightarrow{\expl f^{\dag}}\expl B^{\dag}
\\ &\downarrow &
\\ &\expl X^{\dag}_{P} \end{array}\]
 is a family of curves in $\expl B^{\dag}$ which consist of   stable curves if and only if $f^{\dag}$ is stable.  In particular, exploding a basic log curve of Gross and Siebert or a minimal curve of Abramovich and Chen gives a family of stable exploded curves. To describe what Gross and Siebert's basic condition corresponds to for exploded curves, we shall first need to understand the correspondence between the tropical part  of an exploded manifold and the tropicalization of a  log scheme.
 
 \section{Tropical structure}\label{tropical}

Each exploded manifold $\ex B$ has a tropical part $\totb{\ex B}$, related to the tropicalization of a log scheme $\trop(X^{\dag})$ defined by Gross and Siebert. For example, 
\[\totb{\et mP}=P\]
and
\[\trop{X^{\dag}_{P}}=\hom(Q,[0,\infty))=\bar P\]
where $Q$ is the monoid of integral linear maps $P\longrightarrow [0,\infty)$ and $\bar P$ is the closure of $P\subset\mathbb R^{m}$. A modification of the definition of tropicalization using homomorphisms that are strictly positive on $\alpha^{-1}(0)$ produces $P$ as the tropicalization of $X^{\dag}_{P}$ instead of $\bar P$.

\

As any map $f:X^{\dag}_{P_{1}}\longrightarrow X^{\dag}_{P_{2}}$ involves a homomorphism $Q_{2}\longrightarrow Q_{1}$, it induces an integral linear map
\[\trop(f):\trop X^{\dag}_{P_{1}}\longrightarrow \trop X^{\dag}_{P_{2}}\]
Similarly, if $\et {m_{i}}{P_{i}}$ are represented by log schemes, then any map $f:\et {m_{1}}{P_{1}}\longrightarrow \et {m_{2}}{P_{e}}$ involves a homomorphism $\af {P_{2}}\longrightarrow \af {P_{1}}$ which sends $\afs {P_{2}}$ inside $\afs {P_{1}}$. Such a homomorphism is equivalent to an integral affine map
\[\totb f: P_{1}\longrightarrow P_{2}\] 
 $\trop (f)$ and $\totb{\expl f}$ are the same if $P_{i}$ are closed, and  if $P_{i}$ are not closed, then $\trop(f)$ is obtained from $\totb{\expl f}$ by taking closures.  
 
 \
 
 Suppose that  $(X,\mathcal M_{X})$ is an explodable log scheme. To each point $x\in X$, we may associate the closed integral affine cone $\hom(\bar{\mathcal M}_{X,x},[0,\infty))$, where $\bar {\mathcal M}_{X,x}$ is the stalk of $\bar{\mathcal M}_{X}$ at $x$. Locally $(X,\mathcal M_{X})$ is isomorphic to a  open subset of $X^{\dag} _P$ for which each strata is connected, where each of these cones is either the closure of $P$ or the closure of a face of $P$. We can therefore identify these cones with their image in the closure of $P$, and then define $\trop (X,\mathcal M_{X})$ by taking the quotient of  the disjoint union of these integral affine cones for each $x\in X$ by the equivalence relation generated by these local identifications. 
 
 Similarly, to each point $p\in\totl{\ex B}$ we may associate a polytope $\mathcal P(p)$ so that if an open neighborhood of $p$ is isomorphic to an open subset of $\et mP$, then the image of this point under the tropical part map $\et mP\longrightarrow P$  is in the interior of a face of $P$ isomorphic to $\mathcal P(p)$. The tropical part $\totb{\ex B}$ of $\ex B$ may be defined by taking the quotient of the disjoint union of $\mathcal P(p)$ for each $p\in \totl{\ex B}$ under the equivalence relation generated by locally identifying each $\mathcal P(p)$ with a face of $P$ for each $p$ in a connected open subset of $\ex B$ isomorphic to an open subset of $\et mP$ with connected strata.
 
 \
 
 If $X^{\dag}$ is explodable, then $\trop X^{\dag}$ may be obtained from $\totb{\expl X^{\dag}}$ simply by taking closures of all the polytopes $\mathcal P(p)$ involved. The construction in each case is functorial so given any map 
 \[f:X^{\dag}\longrightarrow Y^{\dag}\] 
there is an induced tropicalization
\[\trop(f):\trop(X^{\dag})\longrightarrow \trop (Y^{\dag})\] 
which is an integral linear map restricted to each of the cones involved.
For exploded manifolds, given any map 
\[g:\ex B\longrightarrow \ex C\]
 $g$ has a tropical part
 \[\totb g:\totb{\ex B}\longrightarrow \totb{\ex C}\]
Note that as a special case of this, each point in $\ex B$ has a tropical part which is a point in the tropical part of $\ex B$, and there is a corresponding surjective map of sets $\ex B\longrightarrow \totb{\ex B}$.

\

For example
\begin{enumerate}
\item
 if $X^{\dag}$ is a toric manifold with log structure given by its toric boundary divisors, then $\trop X^{\dag}=\totb{\expl X^{\dag}}$ may be identified with the toric fan of $X$. 

\

 \item If $X^{\dag}$ is a log scheme associated to a complex manifold with a normal crossing divisor, then $\trop X^{\dag}=\totb{\expl X^{\dag}}$ is the dual intersection complex of $X^{\dag}$, with a single vertex representing the interior of $X^{\dag}$, a ray $[0,\infty)$ representing each irreducible component of the divisor, and a simplex $[0,\infty)^{n}$ representing each $n$-fold intersection of irreducible components.

\

\item If $\ex C$ is a curve, then\begin{itemize}\item if $p\in\totl{\ex C}$ is a node, then $\mathcal P(p)=[0,l]$,\item if $p$ is a marked point of $\totl{\ex C}$, then $\mathcal P(p)=[0,\infty)$,\item and if $p$ is any other point, then $\mathcal P(p)$ is a point. 
\end{itemize}Therefore $\totb{\ex C}$ is a complete integral affine graph where each internal edge corresponds  to a node of $\totl{\ex C}$, each infinite edge  corresponds to a marked point of $\totl{\ex C}$, and each vertex corresponds to a component of the nodal curve $\totl{\ex C}$. Call such a complete integral affine graph a tropical curve.

\

The tropical part of a curve $f:\ex C\longrightarrow \ex B$ in $\ex B$ is an integral affine map 
\[\totb{f}:\totb{\ex C}\longrightarrow \totb{\ex B}\]
Call such  an integral affine map from a tropical curve to $\totb{\ex B}$ a tropical curve in $\totb{\ex B}$.

\

Note that these tropical curves in $\totb{\ex B}$ do not always satisfy the balancing condition of tropical geometry. Instead they satisfy a modified balancing condition determined the underlying map of schemes $\totl f:\totl{\ex C}\longrightarrow \totl{\ex B}$.

\

\item If $(C,\mathcal M_{C})\longrightarrow X^{\dag}_{P}$ is a curve over an explodable log point, then 
\[\begin{array}{c}\trop (C,\mathcal M_{C})
\\ \downarrow
\\ \trop X^{\dag}_{P}:= \bar P\end{array}\]
can be regarded as a family of tropical curves, so the inverse image of any point in $P\subset \bar P$ is a tropical curve. The  length of the internal edge of these tropical curves corresponding to a node $q$ of $C$ is measured by the integral linear function $\rho_{q}$ on $\bar P$. Note that at $0\in \bar P$ and on some faces of $\bar P$, the length of the internal edge corresponding to $q$ shrinks to $0$, but within $P\subset \bar P$, each tropical curve in this family gives an integral affine structure to some fixed graph.

Given a map
\[(C,\mathcal M_{C})\longrightarrow X^{\dag}\]
taking tropicalizations gives a family of tropical curves in $\trop{X^{\dag}}$.
\[\begin{array}{ccc}\trop(C,\mathcal M_{C})&\longrightarrow &\trop X^{\dag}
\\\downarrow
\\ \bar P
\end{array}\]
Of course, if $X^{\dag}$ is explodable,  this family of tropical curves over $P\subset \bar P$ is the tropical part of the family of curves \[\begin{array}{ccc}\expl(C,\mathcal M_{C})&\longrightarrow &\expl X^{\dag}
\\ \downarrow
\\ \expl X^{\dag}_{P}\end{array}\]
\end{enumerate}
 
 \
 
 \
 
In \cite{iec}, an exploded manifold $\ex B$ is called basic if under the above identifications from the definition of $\totb{\ex B}$, $\mathcal P(p)$ is identified with a face of some $\mathcal P(q)$ in at most one way. A stronger condition which is convenient for proving compactness of the moduli space of curves in $\ex B$ is to assume that there exists an integral affine map $\totb{\ex B}\longrightarrow \mathbb R^{N}$ which is injective restricted to each $\mathcal P(p)$. This condition is analogous to the  quasi generated condition from \cite{GSlogGW}. In particular, if $X^{\dag}$ is an explodable log scheme, then $X^{\dag}$ is quasi generated if and only if there exists an integral affine immersion $\totb{\expl X^{\dag}}\longrightarrow \mathbb R^{N}$, and $X^{\dag}$ is almost generated if and only if there exists an integral affine immersion $\totb{\expl X^{\dag}}\longrightarrow [0,\infty)^{N}$.  Using the methods in \cite{GSlogGW} or \cite{cem}, the moduli space of curves  with bounded energy, genus and number of marked points is compact when $X^{\dag}$ is almost generated or $\totb{\ex B}$ is immersible in $[0,\infty)^{N}$, and further combinatorial data is needed for each marked point to get compactness in the case that $X^{\dag}$ is quasi generated, or $\totb{\ex B}$ is immersible in $\mathbb R^{N}$.  The most important example of an exploded manifold with tropical part immersible in $\mathbb R^{N}$ but not $[0,\infty)^{N}$ is $\et m{\mathbb R^{m}}$, which may not be represented by a log scheme.

\

\section{The explosion of basic or minimal curves}\label{basic}

In this section, we shall examine the concepts from \cite{GSlogGW} of the tropical type of a curve and a basic curve. This is related to the minimality condition used by Abramovich and Chen in \cite{acgw}. For simplicity, we shall assume that our log schemes are quasi generated and correspondingly, the tropical parts of our exploded manifolds admit immersions into some $\mathbb R^{N}$.

\

The type defined in \cite{GSlogGW} of a given curve $f:(C,\mathcal M_{C})\longrightarrow (X,\mathcal M_{X})$ records some extra combinatorial information as well as the underlying curve $C\longrightarrow X$. If $(X,\mathcal M_{X})$ is monodromy free in the terminology of \cite{GSlogGW} (in particular if $(X,\mathcal M_{X})$ is quasi generated), then this combinatorial information is  the derivative of the tropical curves in $\trop f$. Explicitly, this information is given as follows:

\

 For each node, $q$, the edge of a tropical curve corresponding to $q$ is equal to $[0,l]$. Taking the derivative of the integral  affine map 
 \[[0,l]\longrightarrow \hom(\bar{\mathcal M}_{X,f(q)},[0,\infty)) \]
  gives a homomorphism
\[u_{q}:\bar{\mathcal M}_{X, f(q)}\longrightarrow \mathbb Z\]
Of course, $u_{q}$ depends on an orientation of this edge - reversing the orientation multiplies $u_{q}$ by $-1$.

\

For each marked point $p$, the edge of a tropical curve corresponding to $p$ is equal to $[0,\infty)$. Taking the derivative of the integral affine map
 \[[0,\infty)\longrightarrow \hom(\bar{\mathcal M}_{X,f(p)},[0,\infty)) \]
gives a homomorphism
\[u_{p}:\bar{\mathcal M}_{X, f(p)}\longrightarrow \mathbb N\]

\

The above description may not  be sufficient to give the tropical type of a curve if $(X,\mathcal M_{X})$ is not monodromy free, because then the map $\hom(\mathcal M_{X,x},[0,\infty))\longrightarrow \trop(X,\mathcal M_{X})$ may involve quotienting by a group action, in which case $\trop f$ may not be sufficient to determine the type of $f$. 

Note that the underlying curve $f:C\longrightarrow X$ is used in the definition of the tropical type to specify the underlying graph for the tropical curves and to specify the stalk of  $\bar{\mathcal M}_{X}$ at the image of nodes and marked points, which chooses out a polytope in $\trop (X,\mathcal M_{X})$ to contain the image of the corresponding edge of the tropical curve.  The following is a definition of tropical type which does not include the information of the underlying curve $f:C\longrightarrow X$.

\begin{defn}\label{tropical type def}
Define a type $\Gamma$ of tropical curve in $\totb{\ex B}$ or $\trop X^{\dag}$ to be a graph $\gamma$ with an orientation of each internal edge and a specification of
\begin{itemize}\item for each vertex $v$ of $\gamma$ a polytope $P_{v}$ in $\totb{\ex B}$ or $\trop X^{\dag}$ 
\item for each edge $e$ of $\gamma$, a polytope $P_{e}$ in $\totb{\ex B}$ or $\trop X^{\dag}$, and  a constant integral vector field $u_{e}$ on $P_{e}$.
\end{itemize}
A tropical curve of  type $\Gamma$ is a tropical curve in $\totb{\ex B}$ or $\trop X^{\dag}$ with an identification of the underlying graph of the tropical curve with $\gamma$ so that the vertex $v$ is sent to the interior of $P_{v}$, the interior of the edge $e$ is sent to the interior of $P_{e}$, and so that if $e$ is an external edge equal to $[0,\infty)$, the map $[0,\infty)\longrightarrow P_{e}$ is an integral curve of $u_{e}$, and if $q$ is an internal edge equal to $[0,l]$ with the positive orientation, the map $[0,l]\longrightarrow P_{q}$ is an integral curve of $u_{q}$.
 
\end{defn}

The relationship of the tropical type of a curve to the marked graph associated to a curve by Chen in \cite{Chen} is as follows: The target space $X^{\dag}$ Chen uses comes with a map $\trop X^{\dag}\longrightarrow [0,\infty)$ which sends every polytope isomorphically to either $0$ or $[0,\infty)$. Then for a node $q$, we may regard $u_{q}$ as giving an integer. If $u_{q}\neq0$, then Chen orients the corresponding edge so that $u_{q}$ is positive and records $u_{q}$ as a contact order $c$. The other edges he leaves unoriented. He also records which vertices are sent to the strata $0$, and doesn't record any information about the external edges $e$ although this is recorded as a contact order $c=u_{e}$ elsewhere.

\

There is a moduli space of tropical curves with a given tropical type. We may parametrize this moduli space by the position of vertices in $p_{v}\in P_{v}$ and the length of internal edges $l_{q}\in(0,\infty)$ . For an internal edge, $q$ of $\gamma$, let 
$v(q,1)$ and $v(q,2)$ be the vertices at the negative and positive boundaries of $e$ respectively.
For there to be any tropical curves of the given type, $P_{v}$ must be a face of $P_{e}$ for each edge $e$ connected to $v$, and if $e$ is an external edge, there must exist an infinite ray in $P_{e}$ in the direction of $u_{e}$. If these conditions are satisfied, 
the condition that there is a tropical curve with data $(p_{v},l_{q})$ is equivalent to the following integral linear equations
\begin{equation}p_{v(q,2)}-p_{v(q,1)}=l_{q}u_{q}\text{ for each internal edge $q$}\end{equation}
For these equations to make sense, identify $P_{v(q,i)}$ as faces of $P_{q}$. Then we may parametrize the moduli space of tropical curves of type $\Gamma$ by the polytope

\[P_{\Gamma}:=\{(p_{v},l_{q})\text{ so that }p_{v(q,2)}-p_{v(q,1)}=l_{q}u_{q} \text{ for all $q$}\}\subset\prod_{v}P_{v}\times\prod_{q}(0,\infty)\]

We can define a universal family of tropical curves  of type $\Gamma$ over $P_{\Gamma}$ with  $(p_{v},l_{q})$ specified by position in $P_{\Gamma}$.

\[\begin{array}{ccc}\hat P_{\Gamma}&\xrightarrow{\underline f_{\Gamma}} &\totb{\ex B}\text{ or }\trop X^{\dag}
\\ \downarrow
\\ P_{\Gamma}\end{array}\]
In other words,  the tropical curve over a point $(p_{v},l_{q})$ in $P_{\Gamma}$ is the unique tropical curve with type $\Gamma$ and the specified $(p_{v},l_{q})$ in $\totb{\ex B}$ or $\trop X^{\dag}$. 

\

The condition of a curve $f$ over a log point in an explodable quasi generated log scheme $X^{\dag}$ being basic is equivalent to the condition that the log point is explodable and the tropical part of the explosion of that curve is a universal tropical curve $\totb f_{\Gamma}$. In other words, $\totb{\expl f}=\totb f_{\Gamma}$, or equivalently $\trop f$ restricted to its interior is equal to $\totb f_{\Gamma}$. If $X^{\dag}$ is Deligne Faltings, then the condition from \cite{acgw} of $f$ being minimal is the same as the condition of $f$ being basic, so it is equivalent to $\totb{\expl f}$ being a universal family of tropical curves.
\[f\text{ basic or minimal }\ \ \ \Longleftrightarrow\ \ \  \totb{\expl f}=\totb f_{\Gamma}\]

\

This condition of being basic or minimal also turns up naturally in the context the moduli stack of curves in an exploded manifold in the following way: Suppose that $P$ is an open polytope, and let 
\[\begin{array}{ccc}\ex C&\xrightarrow{f}&\ex B
\\ \downarrow 
\\ \et mP \end{array}\]
 be a family of curves in $\ex B$ with automorphism group $G$. If $\totb{\ex B}$ is immersible in $\mathbb R^{N}$, then $f/G$ represents a closed substack of the stack of curves in $\ex B$ if and only if $\totb f$ is a universal family $\totb f_{\Gamma}$ of tropical curves in $\totb{ \ex B}$.

\section{(Co)tangent spaces}\label{cotangent}

If $\ex B$ is represented as a log scheme $\ex B^{\dag}$ over $p^{\dag}$, then the cotangent sheaf of $\ex B$ is the same as the sheaf of relative logarithmic differentials $\omega_{\ex B^{\dag}/p^{\dag}}^{1}$ defined in \cite{kk}, so the (co)tangent space of $\ex B$ may be regarded as simply the relative (co)tangent space of $\ex B^{\dag}/p^{\dag}$. In particular, the cotangent sheaf of $\expl X^{\dag}$ is just the pullback of the logarithmic cotangent sheaf $\omega^{1}_{X^{\dag}}$ of $X^{\dag}$.

Explicitly, 
\[\omega_{X^{\dag}_{P}}^{1}=\mathcal O_{X_{P}^{\dag}}\otimes\mathbb Z^{m}\]
where $\mathbb Z^{m}$ is identified with the group of linear integral functions on $P$, which is the group generated by $Q$, the monoid  of nonnegative integral linear functions on $P$.
The sheaf of holomorphic one forms on $\et mP$ is isomorphic to  
\[\mathcal O_{\totl{\et mP}}\otimes \mathbb Z^{m}\]
where $P\subset\mathbb R^{m}$, and $\mathbb Z^{m}$ is identified with integral linear functions on $P$. In terms of the notation in \cite{iec}, the $i$th standard basis vector in $\mathbb Z^{m}$ corresponds to a coordinate function $\tilde z_{i}$ with tropical part the corresponding integral affine function on $P$, and holomorphic differential forms may be locally written as
\[\sum_{i=1}^{m} f_{i}\tilde z_{i}^{-1}d\tilde z_{i}\]
where $f_{i}$ are holomorphic functions.
Similarly, the notation for a holomorphic vectorfield is
\[\sum_{i=1}^{m}f_{i}\tilde z_{i}\frac\partial{\partial \tilde z_{i}}\] 

In terms of exploded manifolds, tangent vectors are defined as derivations. To do this, $\mathcal E^{\times}(\ex B)$ is regarded as the multiplicatively invertible elements in a sheaf $\mathcal E(\ex B)$ of semirings consisting of $\mathbb C\e{\mathbb R}$ valued functions, where the operation of addition is given by 
\[c\e a+d\e b:=\begin{cases}c\e a\text{ if }a<b
\\ (c+d)\e a\text{ if }a=b
\\ d\e b\text{ if }a>b\end{cases}\] 
The tropical part homomorphism $\totb{c\e a}:=a$ is then a homomorphism of semirings onto the tropical semiring, and the smooth part homomorphism 
\[\totl{c\e a}:=\begin{cases}c\text{ if }a=0
\\ 0\text{ if }a>0
\\ \text{ undefined if }a<0\end{cases}\]
is also a homomorphism $\mathbb C\e{[0,\infty)}\longrightarrow \mathbb C$ of semirings. The sheaf of holomorphic maps to $\mathbb C$, (which is $\mathcal O_{\totl{\ex B}}$) may just be constructed from $\mathcal E(\ex B)$ simply by applying the smooth part homomorphism to those functions to which it can be applied.

With the construction of tangent vectors as derivations,  the cotangent complex and any tensor familiar from differential geometry may be constructed.

\

There is one extra notion which makes sense for exploded manifolds but not smooth manifolds: the notion of an integral vector.
\begin{defn}An integral vector in $T_{p}\ex B$ is a vector $v$ at the point $p$  so that for any function $f\in \mathcal E^{\times}$, $vf$ is an integer. Use the notation 
\[{}^{\mathbb Z}T_{p}\ex B\subset T_{p}\ex B\]
for the lattice of integral vectors in $T_{p}\ex B$.

\end{defn}

\

Explicitly, the vector 
\[\sum a_{i}\tilde z_{i}\frac\partial{\partial\tilde z_{i}}\]
at a point $p\in \et mP$ is integral if and only if all the $a_{i}$ are in $\mathbb Z$, and the vector $(a_{1},\dotsc,a_{m})$ at the point $\totb{p}\in P$ is tangent to the interior of some face of $P$.

\

\section{Families}\label{families}

\begin{defn}A family of exploded manifolds is a complete submersion 
\[f:\ex B\longrightarrow \ex F\]
so that $df$ is surjective on integral vectors in the sense that \[df({}^{\mathbb Z}T_{p}\ex B)={}^{\mathbb Z}T_{f(p)}\ex F\]
\end{defn}

\

A family of exploded manifolds is not quite analogous to a proper integral log smooth map. The explosion of an explodable proper log smooth map $f$ is complete if and only if $f$ is integral. The explosion of a proper integral log smooth map is a complete submersion which may send integral vectors from the domain to a proper sub lattice of the integral vectors in the target which has the same $\mathbb Q$-linear span. This is sufficient for a fiber product of an arbitrary map with such a map to always exist, but this fiber product does not behave as nicely as the fiber product of an arbitrary map with a family.

If  $f:(X,\mathcal M_{X})\longrightarrow (Y,\mathcal M_{Y})$ is a map of explodable log schemes, then the condition of $d\expl f$ being surjective on integral vectors is equivalent to the condition that the
 homomorphisms $(\bar{\mathcal M}_{Y,f(x)})^{gr}\longrightarrow (\bar{\mathcal M}_{X,x})^{gr}$ are injective, and have  torsion free cokernel.  Therefore families correspond to proper log smooth maps so that the homomorphisms $(\bar{\mathcal M}_{Y,f(x)})^{gr}\longrightarrow (\bar{\mathcal M}_{X,x})^{gr}$ have torsion free cokernel.

\section{Comparing the moduli stack of curves}\label{moduli stack}

If $C^{\dag}\longrightarrow W^{\dag}$ is a pre-stable family of curves and $W^{\dag}$ is explodable, then $\expl C^{\dag}\longrightarrow \expl W^{\dag}$ is a holomorphic family of curves in the category of exploded manifolds. Any holomorphic family of stable curves in $\expl X^{\dag}$  may be obtained by making a base change on some family $\expl C^{\dag}\longrightarrow \expl W^{\dag}$ obtained using the explosion functor. In this sense, the moduli stack of curves in $\expl X^{\dag}$ is obtained by applying the explosion functor to the moduli stack of curves in $X^{\dag}$. 

\

 For example, if  $\bar {\mathcal M}_{g,n}$ denotes Deligne Mumford space with log structure coming from its boundary divisors, then $\expl\bar{\mathcal M}_{g,n}$ represents the moduli stack of stable exploded curves, and 
\[\begin{array}{c}\expl \bar{\mathcal M}_{g,n+1}\\\downarrow\\ \expl\bar{\mathcal M}_{g,n}\end{array}\] is the universal curve. 

\

The methods used to define a virtual class in \cite{egw} are analytic and use a special type of DeRham cohomology developed in \cite{dre}, however, I expect that in mutually defined cases, Gromov Witten invariants defined in \cite{egw} will agree with those defined in  \cite{acgw} or  \cite{GSlogGW} using the usual construction \cite{BehrendFantechi} with the obstruction theory from \cite{olsson}.  

\

The moduli stack of curves in a family of exploded manifolds is also studied in \cite{cem}, \cite{reg} and \cite{egw}. In \cite{egw}, it is shown that Gromov Witten invariants of exploded manifolds are invariant in connected families and compatible with base change.  As in the absolute case, if the explosion of $X^{\dag}\longrightarrow W^{\dag}$ is a family of exploded manifolds, the moduli stack of curves in $\expl X^{\dag}\longrightarrow \expl W^{\dag}$ may be obtained by applying the explosion functor to the moduli stack of curves in $X^{\dag}/W^{\dag}$.

 An important example is given when $X^{\dag}\longrightarrow X^{\dag}_{[0,\infty)}$ is a simple normal crossings degeneration. The family
 \[\begin{array}{c}\expl X^{\dag}\\\downarrow\\ \et 1{[0,\infty)}\end{array}\]
 of exploded manifolds contains smooth manifolds over points $c\e 0\in\et 1{[0,\infty)}$, and other exploded manifolds which may not be smooth manifolds over  points $c\e x\in \et 1{[0,\infty)}$. The Gromov Witten invariants of any of these exploded manifolds are equal, but tend to be easier to compute in the second situation. All these fibers  over $c\e x$ where $ x>0$ come in one rigid family over $\et 1{(0,\infty)}$ which is the explosion of the log scheme over $(\spec\mathbb  C,\mathbb C^{*}\oplus\mathbb N):=X^{\dag}_{(0,\infty)}$ given by pulling back our family over the obvious map
 \[\begin{array}{ccc}&&X^{\dag}
 \\ &&\downarrow
 \\ X^{\dag}_{(0,\infty)}&\longrightarrow &X^{\dag}_{[0,\infty)}\end{array}\]
 
In some sense, even though $X^{\dag}_{(0,\infty)}$ is called a log `point', it still is a one dimensional object, and the restriction of $X^{\dag}\longrightarrow X^{\dag}_{[0,\infty)}$ to $X^{\dag}_{(0,\infty)}$ is like a family of log schemes, but it does not contain any individual log schemes over $\spec \mathbb C$. Applying the explosion functor to this family gives a one complex dimensional family of exploded manifolds. In the exploded setting, Gromov Witten invariants of a single member of this family may be studied instead of needing to study curves in the whole family. As Gromov-Witten invariants of a single member of this family correspond to Gromov Witten invariants of a smooth fiber, it is most natural in the exploded setting to just study curves in a single member of the family, however it is just as easy to study curves in the entire family, and the same information is obtained.

\section{Reversing the explosion functor}\label{reverse expl}
The following is a list of analogous log objects for exploded manifolds and maps:
\subsection{$\et mP$ for $P$ a cone in $[0,\infty)^{m}$}

\

This is the explosion of  $X^{\dag}_{P}$ defined as in definition \ref{xdef} on page \pageref{xdef}.

\subsection{$\et mP$ for $P\subset [0,\infty)^{n}$ an integral affine polytope}

\

These occur in a family $\et {m+n}{P_{1}}\longrightarrow \et n{P_{2}}$ which is the explosion of some family $X^{\dag}_{P_{1}}\longrightarrow X^{\dag}_{P_{2}}$. Choosing $P_{2}$ to be an open cone allows us to replace $\et mP$ by a log scheme $X^{\dag}_{P_{1}}$ over the log point $(\spec C,\mathbb C^{*}\oplus Q)$ where $Q$ is the monoid of integral linear functions on $P_{2}$.

\subsection{$\et mP$ for $P$ a cone containing an infinite line.}\label{glog}

\

Such a $\et mP$ may be regarded as the explosion of  a generalized log scheme \[X^{\dag}_{P}:=(X,\mathcal M_{X^{\dag}_{P}},\mathcal E^{\times}_{X^{\dag}_{P}})\] which consists of a log scheme $(X,\mathcal M_{X^{\dag}_{P}})$ along with a sheaf of groups $\mathcal E^{\times}_{X^{\dag}_{P}}$ and an inclusion $\mathcal M_{X_{P}^{\dag}}\subset\mathcal E^{\times}_{X_{P}^{\dag}}$ as a sheaf of  saturated sub monoids. With this description $X^{\dag}_{\mathbb R^{n}}$ is given by  
\[ X^{\dag}_{\mathbb R^{n}}:=(\spec \mathbb C,\mathbb C^{*},\mathbb C^{*}\times \mathbb Z^{n})\]
where we regard $\mathbb Z^{n}$ as the integral linear functions on $\mathbb R^{n}$.
Any other integral affine cone $P$ is isomorphic to $P'\times \mathbb R^{n}$ for some $P'\subset [0,\infty)^{m}$. We may characterize $X^{\dag}_{P'\times \mathbb R^{n}}$ as
\[X^{\dag}_{P}=X^{\dag}_{P'}\times X^{\dag}_{\mathbb R^{n}}=(\totl{X^{\dag}_{P'}},\mathcal M_{X^{\dag}_{P'}},(\mathcal M_{X^{\dag}_{P'}})^{gr}\times \mathbb Z^{n})\]

The equivalent of a chart in this case is given by a diagram 
\[\mathbb Z^{m+n}\longleftarrow Q\longrightarrow \mathbb C[Q]/\mathbb C[Q']\]
where $\mathbb Z^{m+n}$ is the group of integral linear functions on $P$, $Q$ is the monoid of nonnegative integral linear functions on $P$ and $Q'$ is the monoid of strictly positive integral linear functions on $P$.

 Of course, when $P\subset[0,\infty)^{n}$, the group of integral linear functions on $P$ is generated by $Q$, so $\mathcal E^{\times}_{X^{\dag}_{P}}=(\mathcal M_{X^{\dag}_{P}})^{gr}$.

We may also regard $X^{\dag}_{P}$ as a stack over the category of explodable log schemes (or another suitably nice category of log schemes). A morphism $X^{\dag}\longrightarrow X^{\dag}_{P}$ (ie an object of this stack $X^{\dag}_{P}$) is equivalent to a homomorphism
from the integral linear functions on $P$ to global sections of $(\mathcal M_{X^{\dag}})^{gr}$ which sends the nonnegative integral linear functions to global sections of $\mathcal M_{X^{\dag}}$ and which sends the strictly positive functions to $\alpha^{-1}(0)\subset\mathcal M_{X^{\dag}}$.  

The easiest case of this is $X^{\dag}_{\mathbb R}$. Objects in this stack are equivalent to log schemes $X^{\dag}$ with a global section of $(\mathcal M_{X^{\dag}})^{gr}$.

\

An example of how such stacks arise is given by the evaluation space $\wedge X^{\dag}$ constructed by Abramovich, Chen, Gillam and Marcus in \cite{acgw}. If $X^{\dag}$ is explodable, then the part of this evaluation space with contact order $c\geq 1$,  $\wedge_{c} X^{\dag}$ is the quotient of some stack $Y^{\dag}$ by the trivial action of the cyclic group of order $c$. This stack $Y^{\dag}$ is locally isomorphic to open subsets of these model stacks $X^{\dag}_{P}$, so the explosion of $Y^{\dag}$ is an exploded manifold. Forgetting this trivial group action and applying the explosion functor to $\wedge X^{\dag}$ gives $\End (\expl X^{\dag})$, one of the evaluation spaces used in \cite{egw}. (The fact that these trivial group actions were thrown out in \cite{egw} leads to some combinatorial correction terms to the gluing formula there.) 

\

A second case in which such generalized explodable log schemes arise is in the gluing formula for Gromov Witten invariants. In particular, the gluing formula for Gromov-Witten invariants of any explodable $X^{\dag}$ with a nontrivial log structure involves Gromov-Witten invariants of these generalized explodable log schemes locally isomorphic to $X^{\dag}_{P}$. This will be discussed further in section \ref{gluing}.

\subsection{A point in an exploded manifold}\label{pt}

\

A point in $\expl X^{\dag}$ will correspond to a point in $X^{\dag}$ if and only if  that point has image at $0\in\totb{\expl X^{\dag}}$. To have a log equivalent, we must always include our point in a family of points $\et m{P}\longrightarrow \expl X^{\dag}$ where $P$ is an open cone, and the tropical part of the map, $P\longrightarrow \totb{\expl X^{\dag}}$ is linear. Then this family of points will be the explosion of a map $X^{\dag}_{P}\longrightarrow X^{\dag}$ from a log point. Of course, there are many choices of $P$, but a canonical choice is to choose the family of all points which have the same smooth part. This will correspond to the log point in $X^{\dag}$ which is gotten by restricting the log structure on $X^{\dag}$ to a point in the underlying scheme.

\subsection{A curve in an exploded manifold}

\

A curve in $\expl X^{\dag}$ will correspond to a curve over $(\spec \mathbb C,\mathbb C^{*})$ in $X^{\dag}$ if and only if its tropical part is a graph with a single vertex and this vertex is sent to $0\in\totb{\expl X^{\dag}}$. Otherwise, we must replace our curve with the moduli space of all curves with the same smooth part and the same tropical type. This will be some family of curves over $\et mP$ for $P$ an open integral affine cone, and this family of exploded curves will be the explosion of a basic or minimal curve over $X^{\dag}_{P}$ in $X^{\dag}$.

\section{Refinements}\label{refinements}

An important operation on exploded manifolds is the operation of refinement. In some sense Gromov-Witten invariants do not change under the operation of refinement. In fact, the virtual moduli space of curves in a refinement is a refinement of the original virtual moduli space of curves.

\begin{defn}A refinement of an exploded manifold $\ex B$ is a complete, bijective submersion 
\[r:\ex B'\longrightarrow \ex B\]
\end{defn}
If the tropical part, $\totb{\ex B}$ of $\ex B$ admits an immersion into $\mathbb R^{N}$, then a refinement of $\ex B$ is equivalent to a subdivision of $\totb{\ex B}$. For example, the explosion of a toric manifold $X^{\dag}$ with its log structure given by its toric boundary divisors is a refinement of $\ex T^{n}_{\mathbb R^{n}}$ given by subdividing $\mathbb R^{n}$ into the toric fan of $X^{\dag}$.

\begin{defn}Let $X^{\dag}$ be a $n$-dimensional toric scheme with its standard log structure coming from its toric boundary divisors. A toric blowup of $X^{\dag}$ is a map 
\[r:Y^{\dag}\longrightarrow X^{\dag}\] 
where $Y^{\dag}$ is a $n$ dimensional toric scheme with its standard log structure, and $r$ is a proper $(\mathbb C^{*})^{n}$-equivariant map which is an isomorphism restricted to the $n$-dimensional toric strata of $Y^{\dag}$ and $X^{\dag}$. 
\end{defn}

\begin{defn} Call a map $r:Y^{\dag}\longrightarrow X^{\dag}$ of log schemes a refinement if it is locally obtained by applying a base change to a toric blowup. \end{defn}

If $r:Y^{\dag}\longrightarrow X^{\dag}$ is a refinement, and $X^{\dag}$ is explodable, then $Y^{\dag}$ is explodable and $\expl r$ is a refinement map. Conversely, if $Y^{\dag}$ and $X^{\dag}$ are explodable, and $\expl r$ is a refinement map, then $r$ is a refinement.

\

 If $Y^{\dag}\longrightarrow X^{\dag}$ is a refinement, then the moduli stack of curves in $Y^{\dag}$ with its obstruction theory should be a refinement of the moduli stack of curves in $X^{\dag}$, with its obstruction theory.  If a cohomology theory which behaves well under refinement is used, then Gromov-Witten invariants of $X^{\dag}$ should be equal to Gromov-Witten invariants of $Y^{\dag}$. I suspect that one case of this is discussed in terms of expanded degenerations by Chen in \cite{Chen2}

\

If $\ex B$ is the explosion of  a complex manifold with normal crossing divisors for which the construction of Ionel in \cite{IonelGW} gives a relative virtual moduli space, I expect that the virtual moduli space of curves in $\ex B$ may be constructed so that the smooth part of a refinement of it is Ionel's moduli space (or at least the closure of the smooth curves in Ionel's moduli space). In the  particular case  that $\ex B$ comes from exploding a complex manifold with a smooth divisor, the smooth part of a refinement of the virtual moduli space of curves in $\ex B$ is equal to the virtual moduli space constructed in \cite{ruan}, \cite{Li}, and \cite{IP}.

\section{Gluing formula}\label{gluing}

The gluing formula for Gromov-Witten invariants defined using log schemes should be equivalent to the gluing formula for Gromov-Witten invariants proved in \cite{egw}.

This gluing formula for Gromov-Witten invariants of $\ex B$ is in the form of a sum over tropical curves in $\totb{\ex B}$, where the contribution of a given tropical curve is determined by taking a fiber product of a virtual moduli space of curves corresponding to each vertex. The rough idea is to construct curves with a given tropical type by gluing together pieces of curves corresponding to each vertex. A difficulty which occurs in the log and exploded settings is that no proper subset of a connected curve is a smooth curve. In \cite{egw}, this difficulty is overcome by using a process called `tropical completion' which takes a piece of a curve which is compact but not complete and completes it to become a smooth curve in a different target, also obtained by the process of tropical completion applied to the closure of the strata of $\ex B$ containing the vertex.

 Whenever a vertex is in the interior of a polytope with dimension $k>0$, the contribution of this vertex involves the virtual moduli space of curves in  an exploded manifold which is some $\et{k}{\mathbb R^{k}}$ bundle, and hence not representable as a log scheme. The corresponding gluing formula for Gromov-Witten invariants of explodable log schemes should include contributions from Gromov Witten invariants of the generalized explodable log schemes defined in section \ref{glog}.

One way to avoid using  generalized log schemes is to refine these generalized log schemes to make them actual log schemes, which involves replacing a $X^{\dag}_{\mathbb R^{k}}$ bundle with a bundle of compact $k$ dimensional toric log schemes. The gluing formula may then use the Gromov-Witten invariants of the corresponding log scheme. As Gromov-Witten invariants should not change under the operation of refinement, it does not matter which $k$ dimensional toric log scheme is used. The gluing theorem from \cite{egw} may then be translated into the language of log schemes.

\

A second way to avoid exploded manifolds which are not log schemes is to work with moduli stacks of compact, but  incomplete curves called cut curves.  I shall attempt to sketch this approach below in a way that translates naturally to the setting of log schemes. Cut curves are a kind of curves with special `cut edges' which can be glued together unambiguously to create usual internal edges or nodes. A cut edge should be regarded as a subset of an internal edge along with the information required for an unambiguous gluing. In the exploded case, this information is given by a special marked point on the cut edge. In the log setting, this special marked point will be translated using the heuristic outlined in section \ref{pt}.  

\

Let $\Gamma$ indicate a  tropical type of curve  in $\totb{\ex B}$ or  $\trop X^{\dag}$ in the sense of definition \ref{tropical type def} from page \pageref{tropical type def}. Let $\Aut \Gamma$ indicate the automorphisms of $\Gamma$. 
Morally, we are interested in the moduli stack of curves with tropical type $\Gamma$, but as this moduli stack is not in general closed, we must define some kind of closure of it.

\begin{defn}A curve almost of type $\Gamma$ is 
\begin{itemize}\item  a basic curve $f$ in $X^{\dag}$ with an  identification of one of the curves in $\trop f$ with a tropical curve of type $\Gamma$,
\item or a curve $f$ in $\ex B$ with an identification as a curve of type $\Gamma$ of one of the curves in the closure of the universal family of tropical curves of type $\totb f$,
\end{itemize}
 \end{defn}

 Use the notation $\mathcal M_{\Gamma}$ for the moduli stack of curves almost of type $\Gamma$. This is morally a closure of the moduli stack of curves with tropical type $\Gamma$, but it may be nonempty even when there are no actual curves of tropical type $\Gamma$. On the other hand, for the virtual moduli space defined as in \cite{egw}, the curves in the virtual moduli space corresponding to $\mathcal M_{\Gamma}$ which do not have tropical type $\Gamma$ have positive codimension.

\

There is an obvious map
\[\mathcal M_{\Gamma}\longrightarrow \mathcal M(\ex B)\text{ or }\mathcal M(X^{\dag})\]
where $\mathcal M(\ex B)$ or $\mathcal M(X^{\dag})$ indicates the moduli stack of curves in $\ex B$ or $X^{\dag}$. When restricted to curves with tropical type $\Gamma$, this map is a degree $\abs{\Aut\Gamma}$ cover of its image. A gluing theorem relates $\mathcal M_{\Gamma}$ to a fiber product of moduli spaces associated with each vertex of $\Gamma$. 

\

In the exploded case, define a moduli stack $\tilde{\mathcal M}_{\Gamma}$ consisting of curves in $\mathcal M_{\Gamma}$ with an extra choice of a point on  each internal edge labeled by an internal edge of $\Gamma$. So $\tilde{\mathcal M}_{\Gamma}\longrightarrow \mathcal M_{\Gamma}$ is a bundle over $\mathcal M_{\Gamma}$ with fiber over a given curve a product of an internal edge for each internal edge of $\Gamma$. Knowledge of $\tilde {\mathcal M}_{\Gamma}$ is equivalent to knowledge of $\mathcal M_{\Gamma}$. We shall describe $\tilde {\mathcal M}_{\Gamma}$ as a fiber product of moduli spaces for each vertex of $\Gamma$. 

Evaluation at one of these extra marked points on an internal edge gives a point in $\ex B$. This defines  an evaluation map
 \[\tilde{\mathcal M}_{\Gamma}\longrightarrow \ex B\]
We can also keep track of the direction of the edge by evaluating the derivative of a curve on the unit integral vector pointing in the positive direction of the edge as oriented by $\Gamma$, which gives an integral vector in ${}^{\mathbb Z}T\ex B$. As $\Gamma$ already keeps track of this vector, we shall not include it in the evaluation map, but simply note here that without this information, the natural target of this evaluation map should be ${}^{\mathbb Z}T\ex B$.

\

 Given a curve in $\tilde{\mathcal M}_{\Gamma}$, the extra information of a point in each internal edge allows us to cut the curve into a different incomplete curve corresponding to each vertex. To do this, we may simply remove all points from the curve with the same image in the tropical part as our chosen point. To remember the information of the chosen point, we must make a definition of a `cut curve' as an abstract exploded space as follows:

 As defined in \cite{iec}, an abstract exploded space is a topological space with a sheaf of groups consisting of $\mathbb C^{*}\e{\mathbb R}$ valued functions which include the constant functions. One way of defining an abstract exploded space is as a subset of an exploded manifold $\ex B$  given the subspace topology and the pullback of the sheaf of functions $\mathcal E^{\times}(\ex B)$. 

\begin{defn} Let  $\et 1{c,(0,l]}$ be the abstract exploded space isomorphic to the subset of $\et 1{(-\infty,l]}$ consisting of the union of the point $1\e 0\in \et 1{(-\infty,l]}$ with all points $c\e a\in\et 1{(-\infty,l]}$ where $a>0$. 
\end{defn}

\begin{defn} A cut curve is a compact abstract exploded space locally isomorphic to an open subset of $\et 1{[0,\infty)}$, $\et 1{[0,l]}$ or $\et 1{c,(0,l]}$. \end{defn}

A cut curve in $\ex B$ is a map of abstract exploded spaces from a cut curve to $\ex B$. Evaluation at a cut is given by evaluating the map at the special point in $\et 1{c,(0,l]}$.

\

An important property of $\et 1{c,(0,l]}$ is that any map of $\et 1{c,(0,l]}$ is entirely determined by its restriction to  $\et 1{(l-\epsilon,l]}\subset \et 1{c,(0,l]}$ for any small enough $\epsilon>0$. In this sense, the deformation theory of a cut curve at a cut edge is similar to the deformation theory of a curve at an external edge or marked point. 

Suppose that a cut tropical curve $\gamma$ in $\totb{\ex B}$ may have its cut edges extended infinitely to give a tropical curve $\gamma'$ in $\totb{\ex B}$.  Then the moduli stack of cut curves with tropical type $\gamma$ is equivalent to the moduli stack of curves with tropical type $\gamma'$ with an extra marked point on each of the external edges of $\gamma'$ which comes from a cut edge.

So long as a cut tropical curve $\Gamma$ has a single vertex, tropical completion may be used to get into the situation above where the cut edges of $\Gamma$ may be extended to be infinite. In particular, $\mathcal M_{\Gamma}$ is a moduli stack of curves in the closure $\ex X$ of some strata of $\ex B$. The tropical completion of $\mathcal M_{\Gamma}$ is the moduli stack of cut curves of type $\Gamma$ in the tropical completion $\check{\ex X}$ of $\ex X$. In this situation, the cut ends of  $\Gamma$ extend infinitely  in $\totb{\check{\ex X}}$, so the moduli space of cut curves of type $\Gamma$ in $\check{\ex X}$ is the same as a moduli space of usual curves with a different tropical type.

\

There are obvious maps 
\[\et 1{c,(0,l_{i}]}\longrightarrow \et 1{[-l_{1},l_{2}]}\]
which realize $\et 1{c,(0,l_{i}]}$ as subsets of $\et 1{[-l_{1},l_{2}]}$ where the special point is $1\e 0$.

\

Given a curve $\ex C$ with a choice of point on an  internal edge, that internal edge is isomorphic to an open subset of  $\et 1{[-l_{1},l_{2}]}$ so that the chosen point is $1\e 0$. Cutting the curve at the given point means replacing the open subset of $\et 1{[-l_{1},l_{2}]}$ with the corresponding open subsets of $\et1 {c,(0,l_{i}]}$ for $i=1,2$.

\

The operation of cutting gives natural maps 
\[\tilde {\mathcal M}_{\Gamma}\longrightarrow \mathcal M_{\Gamma,v}\]
where $\mathcal M_{\Gamma,v}$ indicates  a moduli stack of cut curves in $\ex B$  almost of  tropical type $(\Gamma,v)$ where $(\Gamma,v)$ is determined by cutting $\Gamma$ at internal edges and choosing the connected component containing $v$.
Recall that the internal edges of $\Gamma$ are oriented, so each internal edge is cut into a cut edge with an outgoing orientation and a cut edge with an incoming orientation. Evaluation at outgoing and incoming cut edges gives two different maps
\[\prod_{v}\mathcal M_{\Gamma,v}\longrightarrow \ex B^{k}\]
where $k$ is the number of internal edges of $\Gamma$.
Then $\tilde {\mathcal M}_{\Gamma}$ is the fiber product of $\prod_{v}\mathcal M_{\Gamma,v}$ with itself over the incoming and outgoing evaluation maps. The precise statement for virtual moduli spaces in the exploded case is that the virtual moduli space corresponding to $\tilde{\mathcal M}_{\Gamma}$ is cobordant to the transverse fiber product of the virtual moduli spaces corresponding to $\prod_{v}\mathcal M_{\Gamma,v}$ over the incoming and outgoing evaluation maps. There are obvious equivalents to this gluing theorem which specify the genus and energy of curves at vertices of $\Gamma$.

\

We can now translate this gluing theorem into the language of log schemes.
In the log case, there is no such thing as a point 
on an internal edge, but the edge (node) itself is the equivalent of the moduli space of all points on it.  $\tilde {\mathcal M}_{\Gamma}\longrightarrow\mathcal M_{\Gamma}$ is given as follows: let
\[\begin{array}{ccc}f:&C^{\dag}&\longrightarrow X^{\dag}
\\ &\downarrow
\\ & W^{\dag}\end{array}\]
 be a family in  $\mathcal M_{\Gamma}$, let $C^{\dag}_{q}$ be a node in $C^{\dag}$ labeled by an internal edge $q$ of $\Gamma$. Define $\tilde W^{\dag}$ to be the fiber product of $C^{\dag}_{q}\longrightarrow W^{\dag}$ for each internal edge $q$ of $\Gamma$, and let $\tilde f$ be given by the base change
 \[\begin{array}{ccccc}\tilde f:&\tilde C^{\dag}&\longrightarrow &C^{\dag}&\xrightarrow{f} X^{\dag}
\\ &\downarrow&&\downarrow
\\ &\tilde W^{\dag}&\longrightarrow& W^{\dag}\end{array}\]
Note that the restriction of $f$ to $C^{\dag}_{q}$ gives an evaluation map
\[i_{q}:\tilde W^{\dag}\longrightarrow X^{\dag}\]
 for each internal node $q$, which gives an evaluation map
 \[\tilde {\mathcal M}_{\Gamma}\xrightarrow{i_{q}} X^{\dag}\]

 \

We shall now describe the equivalent of a cut curve.  To get the equivalent of $\et 1{c,(0,l]}$, we must use a different model  consisting of the subset of $\et 1{(-\infty,l]}$ with tropical part $[0,l]\subset(-\infty,l]$ with a special point $1\e 0$. As always, we must  include this in a family containing cut edges of all lengths and all choices of special marked point. This family is the explosion of $X^{\dag}_{c}\longrightarrow X^{\dag}_{(0,\infty)}$ defined below, and the special marked point corresponds to a special section $X^{\dag}_{(0,\infty)}\longrightarrow X^{\dag}_{c}$.

\

\begin{defn} The following is a model
\[\begin{array}{c} X^{\dag}_{c}\\\downarrow \uparrow \\X^{\dag}_{(0,\infty)}\end{array} \]
for a cut edge. 

Define the monoid $M\subset \mathbb N^{2}$ to be all  $(a,b)\in \mathbb N^{2}$ so that if $a=0$, then $b=0$ too. Then let $X^{\dag}_{c}$ be the log scheme with coordinate chart

\[\begin{array}{ccc}M&\longrightarrow&\mathbb C[z]
\\ (a,b)&\mapsto& z^{a}0^{b}\end{array}\]
(with the obvious interpretation of $0^{0}$ as $1$).
Define  the map 
\[X^{\dag}_{c}\longrightarrow X^{\dag}_{(0,\infty)}\] 
in coordinate charts using the diagonal map $\mathbb N\longrightarrow M$, $a\mapsto (a,a)$. Define the special section
\[X^{\dag}_{(0,\infty)}\longrightarrow X^{\dag}_{c}\] 
using the map $M\longrightarrow \mathbb N$ given by $(a,b)\mapsto a$.

\end{defn}

\

As $X^{\dag}_{c}$ is not a fine log scheme, it is useful to have a suitably nice category of log schemes to which it belongs. Below we shall define a category of `$\mathbb R$-saturated' log schemes which contains cut curves and explodable log schemes.  We shall need this category to define families of cut curves.

\begin{defn}
Let $\mathcal M_{x}$ denote the stalk of the structure sheaf of an integral log scheme $X^{\dag}$ at a point $x$. Call a homomorphism 
\[h:\mathcal M^{gr}_{x}\longrightarrow\mathbb R\]
positive if 
\[\alpha^{-1}(0)\subset h^{-1}(0,\infty)\]
and
\[\mathcal M_{x}\subset h^{-1}[0,\infty)\]
In the above, $\alpha^{-1}(0)$ indicates the inverse image of functions which vanish at $x$.
\end{defn}

For example, on $X^{\dag}_{P}$ any point in the closure of a strata of $P$ defines a homomorphism  $\mathcal M_{x}\longrightarrow[0,\infty)$ for any point $x$ in the corresponding strata of the underlying scheme of $X^{\dag}_{P}$. This homomorphism is a positive homomorphism if and only if the point is in the interior of that strata. In this sense, the points of $P$ correspond to positive homomorphisms on $X^{\dag}_{P}$. I think that the tropicalization should be defined using only positive homomorphisms, and not all homomorphisms $\mathcal M_{x}\longrightarrow[0,\infty)$.

\begin{defn}
A $\mathbb R$-saturated log scheme is a saturated integral log scheme $X^{\dag}$ so that for all $x\in X^{\dag}$ an element $f\in\mathcal M_{x}^{gr}$ is in  $\alpha^{-1}(0)\subset\mathcal M_{x}$ if and only if $h(f)>0$ for all positive homomorphisms $h:\mathcal M_{x}^{gr}\longrightarrow \mathbb R$.

\

Say that a saturated integral log scheme  $X^{\dag}$ is $\mathbb R$-saturatable if for all elements $f\in \mathcal M_{x}$ so that $\alpha(f)\neq 0$ there exists at least one positive homomorphism $h:\mathcal M_{x}^{gr}\longrightarrow \mathbb R$ so that $h(f)=0$. We may take the $\mathbb R$-saturation of a $\mathbb R$-saturatable $(X,\mathcal M_{X})$ by adding to $\alpha^{-1}(0)$ all other elements  $g\in\mathcal M_{x}^{gr}$ so that $h(g)>0$ for all positive homomorphisms $h$.
\end{defn}

Note that explodable log schemes, exploded manifolds and $X^{\dag}_{c}$ are $\mathbb R$-saturated log schemes. The following lemma gives the universal property of $\mathbb R$-saturation:

\begin{lemma}\label{universal saturation}Let  $X^{\dag}$ be a $\mathbb R$-saturatable log scheme, and $X^{\dag'}\longrightarrow X^{\dag}$ be its $\mathbb R$-saturation. This has the universal property that if $Y^{\dag}$ is any $\mathbb R$-saturated log scheme, then any map $f:Y^{\dag}\longrightarrow X^{\dag}$ factors uniquely as 
\[Y^{\dag}\xrightarrow{f'} X^{\dag'}\longrightarrow X^{\dag}\]

\end{lemma}
\pf
As $X^{\dag}$ and $X^{\dag'}$ have the same underlying scheme and groupified structure sheaf, we may use the same map on the level of schemes and groupified structure sheaves. All that needs to be verified is that the extra elements of $\alpha^{-1}(0)$ in $\mathcal M_{X^{\dag'}}$ are sent to elements  within $\alpha^{-1}(0)$ of $\mathcal M_{Y^{\dag}}$. This follows from the observation that the pullback of a positive homomorphism $M^{gr}_{y}\longrightarrow \mathbb R$ is a positive homomorphism $\mathcal M^{gr}_{f(y)}\longrightarrow\mathbb R$. 

\stop

\begin{lemma} Let $W^{\dag}$ be a $\mathbb R$-saturated log scheme and the following a fiber product diagram in the category of log schemes.
\[\begin{array}{ccc} C^{\dag}&\longrightarrow & X^{\dag}_{c}
\\ \downarrow &&\downarrow 
\\ W^{\dag}&\longrightarrow &X^{\dag}_{(0,\infty)}\end{array}\] 
Then $C^{\dag}$ is $\mathbb R$-saturatable, and its $\mathbb R$-saturation $C^{\dag'}$ is the fiber product in the category of $\mathbb R$-saturated log schemes
\[\begin{array}{ccc} C^{\dag'}&\longrightarrow & X^{\dag}_{c}
\\ \downarrow &&\downarrow 
\\ W^{\dag}&\longrightarrow &X^{\dag}_{(0,\infty)}\end{array}\]

\end{lemma}

\pf

First, note that the map $X^{\dag}_{c}\longrightarrow X^{\dag}_{(0,\infty)}$ is integral, so $C^{\dag}$ is a saturated integral log scheme.

 At a point $(w,x)$ of the underlying scheme of $C^{\dag}$, $\mathcal M_{(w,x)}$ fits into the following pushout diagram

\[\begin{array}{ccc}\mathcal M_{(w,x)}&\longleftarrow &\mathcal M_{x}
\\ \uparrow && \uparrow
\\ \mathcal M_{w}&\longleftarrow & \mathbb C^{*}\oplus \mathbb N \end{array}\]

If $(f,g)\in\mathcal M_{(w,x)}$ is not sent to $0$ by $\alpha$, then $f\in \mathcal M_{w}$ and $g\in\mathcal M_{x}$ are not sent  to $0$ by $\alpha$, so there exist positive homomorphisms $h_1:\mathcal M_{w}^{gr}\longrightarrow \mathbb R$ and $h_{2}:\mathcal M^{gr}_{x}\longrightarrow \mathbb R$ so that $h_{1}(f)=0=h_{2}(g)$. Each of these positive homomorphisms pull back to a positive homomorphism $\mathbb C^{*}\oplus \mathbb Z\longrightarrow \mathbb R$, which is determined by a positive number which is the image of $(1,1)$. We may therefore multiply $h_{1}$ by a positive constant to obtain a positive homomorphism which vanishes on $f$ and which restricts to be the same as $h_{2}$ on $\mathbb C^{*}\oplus\mathbb Z$. This defines a positive homomorphism $\mathcal M_{(w,x)}^{gr}\longrightarrow \mathbb R$ which vanishes on $(f,g)$. It follows that $C^{\dag}$ is $\mathbb R$-saturatable. 

Lemma \ref{universal saturation} implies that the $\mathbb R$-saturation $C^{\dag'}$ of $C^{\dag}$ has the required universal property to be the fiber product in the category of $\mathbb R$-saturated log schemes.

\stop

It is important to note that when a (fiber) product is used below, it means the (fiber) product in the category of $\mathbb R$-saturated log schemes, which in general is the usual (fiber) product followed by the operation of $\mathbb R$-saturation.

\begin{defn} A cut (log) curve is a proper map $C^{\dag}\longrightarrow W^{\dag}$ in the category of $\mathbb R$-saturated log schemes  with a special section $W^{\dag}\longrightarrow C^{\dag}$ for each cut edge, so that $C^{\dag}\longrightarrow W^{\dag}$ with its special sections is
 locally isomorphic to an open subset of the base change (in the category of $\mathbb R$-saturated log schemes) of one of the following models:
\begin{enumerate}
\item A cut edge \[\begin{array}{c}X^{\dag}_{c}\\ \downarrow\uparrow\\ X^{\dag}_{(0,\infty)} \end{array}\]
\item A node \[\begin{array}{c}X^{\dag}_{[0,\infty)^{2}}\\ \downarrow\\X^{\dag}_{[0,\infty)} \end{array}\]
\item A marked point
\[\begin{array}{c}X^{\dag}_{[0,\infty)}\\ \downarrow\\ (\spec \mathbb C, \mathbb C^{*})\end{array}\]
\end{enumerate}
A cut curve in $X^{\dag}$ is a map $f:C^{\dag}\longrightarrow X^{\dag}$.

 Evaluation at a cut edge is given by the composition of $f$ with the special section $W^{\dag}\longrightarrow C^{\dag}$ corresponding to the pullback of the special section $X^{\dag}_{(0,\infty)}\longrightarrow X^{\dag}_{c}$ at the cut edge.

\end{defn}

As the explosion of the moduli stack of cut curves in $X^{\dag}$ is the moduli stack of cut curves in $\expl X^{\dag}$ which has a nice obstruction theory allowing the definition of a virtual fundamental class,  the moduli stack of cut curves in $X^{\dag}$ should have nice obstruction theory allowing the definition of a virtual fundamental class.   Let $(\Gamma,v)$ indicate the type of cut tropical curve given by cutting $\Gamma$ at internal edges and choosing the component containing $v$. Denote by $\mathcal M_{\Gamma,v}$ the moduli stack of cut curves in $X^{\dag}$ which almost have tropical type $(\Gamma,v)$.

Evaluation at incoming and outgoing cut edges gives two maps 
\[\prod _{v}\mathcal M_{\Gamma,v}\longrightarrow (X^{\dag})^{k}\]
where $k$ is the number of internal edges of $\Gamma$.

Then $\tilde {\mathcal M}_{\Gamma}$ should be a fiber product of $\prod_{v}\mathcal M_{\Gamma,v}$ with itself over the outgoing and incoming evaluation maps, and the same should hold for virtual fundamental classes. 

The virtual fundamental class of $\mathcal M_{\Gamma,v}$ should only depend on the closure of the strata of $X^{\dag}$ corresponding to $v$. This closure of a strata will be equal to a $X^{\dag}_{P}$ bundle over some log smooth scheme, where $P$ is an open cone in some $\mathbb R^{m}$. Actually, the virtual fundamental class should only depend on the corresponding  $X^{\dag}_{\mathbb R^{m}}$ bundle, which corresponds to the tropical completion used in \cite{egw}.

\

\

The cutting map $\tilde {\mathcal M}_{\Gamma}\longrightarrow \prod_{v}\mathcal M_{\Gamma,v}$ may be described as follows:

\

 A neighborhood of a node in $\tilde C^{\dag}$ is a base change of the map on the left hand side of the following pullback diagram 
 \[\begin{array}{ccc}X^{\dag}_{P}&\longrightarrow &X^{\dag}_{[0,\infty)^{2}}
 \\\downarrow &&\downarrow
 \\  X^{\dag}_{(0,\infty)^{2}}&\longrightarrow & X^{\dag}_{[0,\infty)}\end{array}\]
where the tropical parts of the bottom and right hand map are $x_{1}+x_{2}$. The polytope $P$ is  a cone over $[0,1]\times (0,1)$. The diagonal section $X^{\dag}_{(0,\infty)^{2}}\longrightarrow X^{\dag}_{P}$ corresponding to the obvious inclusion $X^{\dag}_{(0,\infty)}\subset X^{\dag}_{[0,\infty)}$ pulls back to a section of $\tilde C^{\dag}\longrightarrow \tilde W^{\dag}$ which corresponds in the exploded setting to the section of an edge corresponding to an extra marked point. 

There are two inclusions of $X^{\dag}_{c}\times X^{\dag}_{(0,\infty)}$ into $X^{\dag}_{P}$ with tropical image the subset of $P$ above or below the diagonal. We shall describe the inclusion into $X^{\dag}_{P}$ with tropical image the subset of $P$ below the diagonal, the other inclusion being analogous.

\

Describe $X^{\dag}_{P}$ using the chart
\[(a,b,c,d)\mapsto x^{a}y^{b}0^{c}0^{d}\in \mathbb C[x,y]/\{xy=0\}\]
where $(a,b,c,d)$ is in the quotient of  $\mathbb N^{4}$ by the equivalence relation  generated by $(1,1,0,0)=(0,0,1,1)$. The map $X^{\dag}_{P}\longrightarrow X^{\dag}_{(0,\infty)^{2}}$ is given by \[(c,d)\mapsto (0,0,c,d)\] and the section $X^{\dag}_{(0,\infty)^{2}}\longrightarrow X^{\dag}_{P}$ is given by \[(a,b,c,d)\mapsto (a+c,b+d)\]
 
 \
 
 Describe $X^{\dag}_{c}\times X^{\dag}_{(0,\infty)}$ by the chart
\[(a,b,c)\mapsto x^{a}0^{b+c}\]
where $(a,b,c)$ are in  the submonoid  of $\mathbb N^{3}$ so that if $a=0=c$, then $b=0$ too. Note that this monoid is not $M\times \mathbb N$ as we are taking the product in the category of $\mathbb R$-saturated log shemes and not the category of arbitrary log schemes! The map $X^{\dag}_{c}\times X^{\dag}_{(0,\infty)}\longrightarrow X^{\dag}_{(0,\infty)^{2}}$ is given by \[(a,b)\mapsto (a,a,b)\] The canonical section is given by \[(a,b,c)\mapsto (a,c)\]

\

The inclusion $X^{\dag}_{c}\times X^{\dag}_{(0,\infty)}\longrightarrow X^{\dag}_{P}$ is given in the above charts by the map 
\begin{equation}\label{abcd}(a,b,c,d)\mapsto (a+c,b+c,b+d)\end{equation} 

Verification that this inclusion commutes with the maps to $X^{\dag}_{(0,\infty)^{2}}$ and sections is immediate.

\begin{lemma}\label{cutsubset} $X^{\dag}_{c}\times X^{\dag}_{(0,\infty)}$ may be regarded as the subset of $X^{\dag}_{P}$ with tropical image the diagonal and every thing below it in $P$ in the sense that any map of a $\mathbb R$-saturated log scheme to $X^{\dag}_{P}$ with tropical image contained in this subset of $P$ factors through the inclusion $X^{\dag}_{c}\times X^{\dag}_{(0,\infty)}\longrightarrow X^{\dag}_{P}$.  

\end{lemma}

\pf

The tropicalization of $X^{\dag}_{P}$ is the closure of $P$, which corresponds to homomorphisms $\mathbb N^{4}\longrightarrow[0,\infty)$ so that $(1,1,0,0)$ is given the same value as $(0,0,1,1)$. These homomorphisms may be written (nonuniquely) in the form $(a,b,c,d)\mapsto x_{1}(a+c)+x_{2}(b+c)+x_{3}(b+d)+x_{4}(a+d)$ where $x_{i}\geq 0$. The homomorphisms in (the closure of) our subset of $P$ are those which are non negative on $(-1,0,1,0)$, which are those that may be written as $x_{1}(a+c)+x_{2}(b+c)+x_{3}(b+d)$ where $x_{i}\geq 0$.

The  homomorphisms $\mathcal M_{X^{\dag}_{P}}^{gr}\longrightarrow \mathbb R$ which restrict to be positive homomorphisms somewhere correspond to the points of $P$. Our subset of $P$ corresponds to a subset of positive homomorphisms $S$.  In our coordinate chart, these positive homomorphisms are precisely the homomorphisms  in the form
\[(a,b,c,d)\mapsto x_{1}(a+c)+x_{2}(b+c)+x_{3}(b+d)\]
 where $x_{1}\geq 0$, $x_{2}\geq 0$, $x_{1}+x_{2}>0$ and $x_{3}>0$. 
These are the positive homomorphisms from $X^{\dag}_{P}$ which have the extra property  of being non negative on $(-1,0,1,0)$.  Note that these are precisely the pull back of the positive homomorphisms from $X^{\dag}_{c}$ under the map (\ref{abcd}).

A map $Y^{\dag}\longrightarrow X^{\dag}_{P}$ has tropical image contained in the required subset if and only if for all points $y$ in the underlying scheme of $Y^{\dag}$, the pullback of any homomorphism $\mathcal M_{y}\longrightarrow [0,\infty)$ corresponds to a homomorphism which is nonnegative on $(-1,0,1,0)$. Therefore, any positive homomorphism $\mathcal M^{gr}_{y}\longrightarrow \mathbb R$ pulls back to one of these positive homomorphisms in $S$. 

Any such map must have image in the component of the underlying scheme of $X^{\dag}_{P}$ which is the image of $X^{\dag}_{c}$, so this underlying map of schemes factors through the underlying map of schemes in $X^{\dag}_{c}\times X^{\dag}_{(0,\infty)}\longrightarrow X^{\dag}_{P}$.

For any point $x$ in this underlying scheme of $X^{\dag}_{c}\times X^{\dag}_{(0,\infty)}$, let  $\mathcal M_{x}$ indicate the stalk of the defining sheaf of $X^{\dag}_{P}$ at $x$ and $\mathcal M_{x}'$ denote the stalk of the defining sheaf of $X^{\dag}_{c}\times X^{\dag}_{(0,\infty)}$ at this point. We have an inclusion $\mathcal M_{x}\longrightarrow \mathcal M_{x}'$ which extends to an isomorphism of groups. The extra elements of $M_{x}'\subset\mathcal M_{x}^{gr}$ are those elements on which any positive homomorphism from $S$ is strictly positive. As $Y^{\dag}$ is $\mathbb R$-saturated, these elements are sent to elements of $\mathcal M_{y}$ which are sent to $0$ by $\alpha$. 

It follows that there is a unique map $Y^{\dag}\longrightarrow X^{\dag}_{c}\times X^{\dag}_{{(0,\infty)}}$ which is the same as the original map on the level of underlying schemes and groupification of structure sheaves.

This  factorizes  our original map into 
\[Y^{\dag}\longrightarrow X^{\dag}_{c}\times X^{\dag}_{(0,\infty)
}\longrightarrow X^{\dag}_{P}\]
as required.

\stop

\

We may similarly define $X^{\dag}_{c}\times X^{\dag}_{(0,\infty)}$ as the subset of $X^{\dag}_{P}$ with tropical image everything above the diagonal. 

\

The cutting construction is now easy to describe. After base change, a neighborhood of a node with its tautological section is locally isomorphic to an open subset of a base change of  $X^{\dag}_{P}\longrightarrow X^{\dag}_{(0,\infty)^{2}}$. The two subsets of $X^{\dag}_{P}$ corresponding to the two inclusions of $X^{\dag}_{c}\times X^{\dag}_{(0,\infty)}$ give two subsets of the corresponding base change. Cutting replaces the original chart with these two subsets. Lemma \ref{cutsubset} implies that this construction is functorial, so the construction may be globalized.

\section{Analogue of log schemes in the $C^{\infty}$ setting}\label{smooth log}

So far, this paper has compared log schemes with exploded manifolds in the case of integrable complex structures. The reader who likes both log schemes and analysis might be wondering what would be analogous to a log scheme in the $C^{\infty}$ setting.  

To view a complex manifold  $X$ as a $C^{\infty}$ manifold, the sheaf $\mathcal O_{X}$ of holomorphic functions should be replaced by a sheaf of smooth functions. It is convenient to have an inclusion of $\mathcal O_{X}$ into this sheaf of smooth functions, so $\mathcal O_{X}$ should be replaced by a sheaf $\mathcal C_{X}$ of smooth maps to $\mathbb C$. We may define $\mathcal C_{X}$ for any complex variety as follows:

\begin{defn}\label{CX}If $X$ is a complex variety, define the sheaf $\mathcal C_{X}$ of $C^{\infty}$ maps to $\mathbb C$ to be the sheaf of $\mathbb C$ valued functions locally equal to 
\[g(\zeta_{1},\dotsc,\zeta_{n})\]
where $\zeta_{i}$ are locally defined functions in $\mathcal O_{X}$, and $g:\mathbb C^{n}\longrightarrow \mathbb C$ is $C^{\infty}$.
\end{defn}

Of course, if a variety $X$ is viewed as locally given by a subset of $\mathbb C^{n}$, these functions in $\mathcal C_{X}$ are the restriction of smooth complex valued functions on $\mathbb C^{n}$ to $X$. A similar definition of $\mathcal C_{X}$ for $X$ a scheme over $\mathbb C$ should be possible but as elements of $\mathcal O_{X}$ are not correctly viewed as functions on the set of points in $X$, more care must be taken. The underlying schemes of $X^{\dag}_{P}$ are always varieties, so the above construction suffices to construct the sheaf of smooth complex valued functions on the underlying scheme of $X^{\dag}_{P}$. 

If $(X,\mathcal M_{X})$ is a variety with a log structure, we may construct another sheaf of monoids $\mathcal N_{X}$ using the following pushout diagram

\[\begin{array}{ccc}\mathcal N_{X}&\longleftarrow &\mathcal M_{X}
\\\uparrow &&\uparrow
\\ \mathcal C^{*}_{X}&\longleftarrow&\mathcal O^{*}_{X}\end{array}\] 
In other words, elements of $\mathcal N_{X}$ are in the form of $fg$ for $f\in\mathcal C^{*}_{X}$ and $g\in\mathcal M_{X}$, and the obvious equivalence of $(fh)g=f(hg)$ holds if $h\in O^{*}_{X}$. Using the map $\alpha:\mathcal M_{X}\longrightarrow \mathcal O_{X}\subset\mathcal C_{X}$, we can define a map $\mathcal N_{X}\longrightarrow \mathcal C_{X}$ by $fg\mapsto f\alpha (g)$.

The topological space $X$ along with $\mathcal N_{X}\longrightarrow \mathcal C_{X}$ should be thought of as giving the $C^{\infty}$ analogue of $(X,\mathcal M_{X})$.  We can make the following rather general definition of a $C^{\infty}$ log space:

\begin{defn} A $C^{\infty}$ log space $X$ is a topological space $X$ with
\begin{enumerate}
\item A sheaf $\mathcal C_{X}$ of maps to $\mathbb C$ so that
\begin{enumerate}
\item if $g:\mathbb C^{n}\longrightarrow \mathbb C$ is smooth and $\zeta_{1},\dotsc,\zeta_{n}$ are in $\mathcal C_{X}(U)$, then $g(\zeta_{1},\dotsc,\zeta_{n})$ is in $\mathcal C_{X}(U)$.
\item The topology on $X$ is given by the functions in $\mathcal C_{X}$ in the sense that every function in $\mathcal C_{X}$ is continuous and every closed subset of $X$ is the zero set of some globally defined function in $\mathcal C_{X}$.
\end{enumerate}
\item A sheaf of monoids $\mathcal N_{X}$ on $X$.
\item A map of sheaves of monoids 
\[\alpha:\mathcal N_{X}\longrightarrow (\mathcal C_{X},\times)\]
which is an isomorphism when restricted to the inverse image of the sheaf $\mathcal C^{*}_{X}$ of functions in $\mathcal C_{X}$ which are nowhere zero.
\end{enumerate}

\end{defn}

The notion of a $C^{\infty}$ space as a topological space  $X$ with a sheaf $\mathcal C_{X}$ of smooth functions should be replaced by a more sophisticated notion if a $C^{\infty}$ analogue of a scheme is desired, or replaced by a stronger notion if a `nice' space is wanted. Beyond this choice of what to call a $C^{\infty}$ space, the above definition is a straightforward copy of the definition of a log scheme using the sheaf of smooth $\mathbb C$ valued functions instead of $\mathcal O_{X}$. If $\mathbb R$ is used instead of $\mathbb C$, then it is natural to put a log structure on manifolds with boundary and corners as was noted in  \cite{gokova}.

\begin{defn}
Suppose that $(X,\mathcal M_{X})$ is a complex variety with a log structure. Then define the $C^{\infty}$ analogue  
\[ \text{Smooth} (X,\mathcal M_{X})\]
of $(X,\mathcal M_{X})$ as follows: Smooth$(X,\mathcal M_{X})$ is the $C^{\infty}$ log space which is the topological space $X$ with
\begin{itemize}\item the sheaf $\mathcal C_{X}$ of smooth maps to $\mathbb C$ from Definition \ref{CX} 
\item  a sheaf $\mathcal N_{X}$ of monoids constructed via the following pushout diagram
\[\begin{array}{ccc}\mathcal N_{X}&\longleftarrow &\mathcal M_{X}
\\\uparrow &&\uparrow
\\ \mathcal C^{*}_{X}&\longleftarrow&\mathcal O^{*}_{X}\end{array}\] 
\item the map 
\[\mathcal N_{X}\longrightarrow \mathcal C_{X}\]
induced from the map $\alpha:\mathcal M_{X}\longrightarrow \mathcal O_{X}\subset \mathcal C_{X}$
and the inclusion $\mathcal C^{*}_{X}\subset\mathcal C_{X}$.
\end{itemize}
\end{defn}

Nice $C^{\infty}$ log spaces are given by Smooth$(X^{\dag}_{P})$ where $X^{\dag}_{P}$ is as in definition \ref{xdef}. Other examples of well behaved $C^{\infty}$ log spaces are smooth manifolds $M$ with $\mathcal C_{M}$ the space of smooth $\mathbb C$  valued functions and $\mathcal N_{M}:=\mathcal C^{*}_{M}$.

\begin{defn}A $C^{\infty}$ log space $(X,\mathcal N_{X}\rightarrow \mathcal C_{X})$ is explodable if it is locally isomorphic to an open subset of a product of $\mathbb R^{n}$ with Smooth$(X^{\dag}_{P})$.

\end{defn}

Important examples of such explodable $C^{\infty}$ log spaces will be constructed in the next section to correspond to almost complex manifolds with `normal crossing divisors'.

Note that the special log point $p^{\dag}$ is a $C^{\infty}$ log space with $\mathcal C_{p^{\dag}}=\mathbb C$, $\mathcal N_{p^{\dag}}=\mathbb C^{*}\e{[0,\infty)}$ and the map $\mathcal N_{p^{\dag}}\longrightarrow \mathcal C_{p^{\dag}}$ just the smooth part homomorphism defined on page \pageref{smooth part def}. As with explodable log schemes, the explosion of an explodable $C^{\infty}$ log space is  given by the base change

\[\begin{array}{ccc}\expl(X,\mathcal N_{X}\rightarrow \mathcal C_{X})&\longrightarrow &(X,\mathcal N_{X}\rightarrow \mathcal C_{X})
\\ \downarrow&&\downarrow
\\ p^{\dag}&\longrightarrow &p \end{array}\] 
where $p$ indicates a point considered as a $C^{\infty}$ log space.

With the holomorphic case from section \ref{expl} understood, it is straightforward to verify that $\expl(\mathbb R^{n}\times \text{Smooth}X^{\dag}_{P})$ is isomorphic to $\mathbb R^{n}\times \et mP$ viewed as a smooth exploded manifold as defined in \cite{iec}. It follows that the explosion of a explodable $C^{\infty}$ log space is a smooth exploded manifold. 

\subsection{An almost complex analogue of a normal crossing divisor}\label{dbarlog ncd}

\

The cotangent sheaf of a log scheme $(X,\mathcal M_{X})$ has a notion of  $d\log f$ for $f$ a section of $\mathcal M_{X}$. This has the property that $\alpha (f)d\log f=d\alpha(f)$. Any appropriate definition of the sheaf of $\mathbb C$-valued one forms  on a $C^{\infty}$ log space should allow the definition of $d\log f$ for any section of $\mathcal M_{X}$. I expect that the program of Gross and Siebert and Abramovich and Chen may be imitated for holomorphic curves in explodable $C^{\infty}$ log spaces $(X,\mathcal N_{X}\rightarrow \mathcal C_{X})$ that have an almost complex structure $J$ so that $d\log f+d\log (if)\circ J$ is a smooth complex valued one form on $X$ for all sections $f$ of $\mathcal N_{X}$.  In particular, if $g$ is a smooth $\mathbb C$-valued function on $X$ which is in the image of $\alpha:\mathcal N_{X}\longrightarrow\mathcal C_{X}$, then $\dbar g$ should be a smooth complex valued one form times $g$.

\

 \begin{defn} A $\dbar\log$ compatible function on an almost complex manifold $(M,J)$ is a smooth map $f:M\longrightarrow \mathbb C$ so that there exists some smooth complex valued one form $\theta$ so that 
 
 \[\frac 12(df+idf\circ J)=f\theta\]
 Define $\dbar\log f$ to be $\theta$ if $f$ is not identically $0$.
 
 \end{defn}
 
 Corollary \ref{isolated zero} below shall imply  that $\dbar\log f$ is well defined if $f$ is $\dbar\log$ compatible and not identically $0$. 
 
 Any smooth $\mathbb C^{*}$ valued function is $\dbar\log$ compatible.  If the dimension of $M$ is larger than $2$, then  $\dbar\log$ compatible functions which are somewhere $0$   may have to be identically $0$ for generic almost complex structures, however the almost complex structures that do admit nontrivial $\dbar\log$ compatible functions with zeros are sufficiently flexible to define divisors in the almost complex setting as the zero sets of $\dbar\log$ compatible functions.

 \begin{defn}\label{ncd}A $\dbar\log$-compatible normal crossing divisor in an almost complex manifold $(M^{2n},J)$ is a subset $V\subset M$ so that around any point $p\in M$, there exists a neighborhood $U$ and $\dbar\log$ compatible functions 
 \[z_{1},\dotsc,z_{n}:U\longrightarrow \mathbb C\]
 so that 
 \begin{itemize}
 \item $V\cap U=\{z_{1}z_{2}\dotsb z_{n}=0\}$
\item $0$ is a regular value of $z_{i}$
 \item The intersection of any number of the submanifolds $\{z_{i}=0\}$ is transverse.
		  
  \end{itemize}
  A $\dbar\log$-compatible simple normal crossing divisor is a $\dbar\log$-compatible normal crossing divisor $V\subset M$ which is equal a finite union of  proper submanifolds of $M$.
  \end{defn}
 
 \begin{remark}Note that $V$ is a real codimension $2$ holomorphic, proper, immersed submanifold of $M$ which intersects  itself transversely.
 \end{remark}
 
 We shall show in a series of lemmas below that the sheaf $\mathcal N(M,V)$ of $\dbar\log$ compatible functions which are non-vanishing away from $V$ gives $(M,V)$ the structure of an explodable $C^{\infty}$ log space, so if $J$ is tamed by a symplectic form and $V$ is simple, \cite{egw} may be used to define Gromov-Witten invariants of $M$ relative to $V$ (or alternatively, the brave reader could try to construct Gromov-Witten invariants relative to $V$ in analogy to log Gromov-Witten invariants.) 
 
 The notion of a $\dbar\log$ compatible normal crossing divisor is stronger than the definition of a normal crossing divisor used by Ionel in \cite{IonelGW}, but it is still  flexible enough to be able to define Gromov-Witten invariants relative to an appropriate notion of a `symplectic normal crossing divisor', because the obstruction to making an immersed codimension $2$ submanifold a $\dbar\log$ compatible divisor for some tamed almost complex structure is the same as the obstruction to making it a holomorphic submanifold for some tamed almost complex structure. In particular, 
 if $V$ is a finite union of transversely intersecting, proper real codimension $2$ submanifolds of a manifold $M$ with symplectic structure $\omega$ and $J_{t}$ is any finite dimensional family of almost complex structures tamed by $\omega$ so that $V$ is $J_{t}$ holomorphic, then there exists a family $J^{'}_{t}$ of almost complex structures tamed by $\omega$ for which $V$ is a $\dbar\log$ compatible normal crossing divisor and $J'_{t}=J$ restricted to $T_{p}M$ for all $p$ in $V$. Moreover, any choice of  $J'_{t}$ on a compact subfamily extends to a smooth choice on the whole family. In the special case that the symplectic submanifolds making up $V$  intersect symplectically orthogonally, the space of almost complex structures tamed by $\omega$ for which $V$ is holomorphic is nonempty and contractible, however it is not clear to me that the same is true for more general $V$.

  \

  We shall now proceed to prove that a $\dbar\log$ compatible normal crossing divisor induces a $C^{\infty}$ log structure.
 
\begin{lemma}\label{dbarlog division} Suppose that $f$ and $z$ are $\dbar\log$ compatible functions so that $0$ is a regular value of $z$ and $z^{-1}(0)\subset f^{-1}(0)$. Then there exists a $\dbar\log$ compatible function $g$ so that $f=gz$.
 \end{lemma}
 \pf
Let $g$ be equal to $f/z$ where $z\neq 0$. The fact that $0$ is a regular value of $z$ and $z^{-1}(0)\subset f^{-1}(0)$ implies that $g$ is locally bounded and locally gives a well defined element of $L^{p}$ for any $p$. Away from $\{z=0\}$, $g$ satisfies the equation
\[\dbar g=g(\dbar\log f-\dbar\log z)\]
therefore, away from $z^{-1}(0)$,
\[d*\lrb{\lrb{\dbar -(\dbar\log f-\dbar \log z)}g}=0\] 
 where $*$ indicates the hodge star operator using some metric. The operator $\Delta:=d*(\dbar -\dbar\log f +\dbar\log z)$ is elliptic, and $g$ is a weak solution to $\Delta g=0$, so elliptic regularity for $\Delta$ implies that $g$ must be smooth. Therefore $g$ is $\dbar\log$ compatible with $\dbar\log g=\dbar\log f-\dbar\log z$, and $f=gz$.

 \stop
 
 \begin{lemma}\label{zero extension} If $f$ is a $\dbar\log$ compatible function on the complex unit disk so that $f$ and all derivatives vanish at $0$, then $f$ is identically $0$.
 \end{lemma}

 \pf
 
 Let $f$ be such a $\dbar\log$ compatible function which vanishes to infinite order at $0$ so that $\dbar f=gf$.
 Lemma \ref{dbarlog division} implies that $f=f_{1}z$ where $\dbar f_{1}=g f_{1}$. As all derivatives of $f$ vanish at $0$, $f_{1}$ and all derivatives of $f_{1}$ vanish at $0$. Therefore, for all $n$, 
 \[f=f_{n}z^{n}\]
 where 
 \[\dbar f_{n}= f_{n}g\]
 
 Let $B_{r}$ indicate the disk of radius $r$ around $0$. Cauchy's integral formula gives for any  $c\in B_{r}$,
 
	 \[f_{n}(c)=\frac 1{2\pi i}\lrb{\int_{\partial B_{r}} \frac {f_{n}}{z-c}dz+\int_{B_{r}} \frac{f_{n}}{z-c}d z \wedge g}\]
 
 Because $ {\abs z}^{-1}$ is integrable on the complex plane, we may choose $r$ small enough so that the integral of $\abs z^{-1}$ on any ball of radius $r$ is  less than $\pi$ divided by the supremum of $\abs g$ in $B_{r}$. Then for any smooth function $h$,
 \[\frac 1{2\pi }\abs{\int_{B_{r}}\frac h{z-c}dz\wedge g}\leq \frac 12\max_{B_{r}} \abs h\]
Then  
\[\abs {f_{n}(c)}\leq \max_{\partial B_{r}}r(r-\abs c)^{-1}\abs {f_{n}}+\frac 12 \max_{B_{r}}\abs {f_{n}}\]	 
	so
	
	\[\max_{B_{r/2}} \abs{f_{n}}\leq 2\max_{\partial B_{r} }\abs{f_{n}}+\frac 12\max_{B_{r}}\abs {f_{n}}\]
	As $f_{n}=fz^{-n}$, either $f$ is identically $0$ in $B_{r/2}$ or for $n$ large enough the maximum of $\abs{f_{n}}$ in $B_{r}$ is achieved within $B_{r/2}$. Therefore, for $n$ large enough,
	\[\max_{B_{r/2}} \abs{f_{n}}\leq 4\max_{\partial B_{r} }\abs{f_{n}}\]
	 so 
	 \[\max_{B_{r/2}}\abs f\leq \frac 4{2^{n}}\max_{\partial B_{r}}\abs f\]
	 As this is true for all $n$ large enough, $f$ vanishes identically inside $B_{r/2}$. It follows that the set where $f$ and all derivatives vanish is both open and closed.
	 
	 \stop
	 
\begin{cor}\label{isolated zero} If $f$ is a $\dbar\log$ compatible function on the unit complex disk, then the zeros of $f$ are isolated if $f$ is not identically $0$. \end{cor}
	 
	 \pf
	 
	 Suppose without losing generality that the zeros of $f$ have a limit point at $0$. Then $f$ vanishes at zero so $fz^{-1}$ is smooth and $\dbar\log$ compatible and also has $0$ as a limit point of its zeros. It follows that $fz^{-n}$ is smooth for all $n$, so $f$ and all its derivatives vanish at zero. Therefore Lemma \ref{zero extension} implies that $f$ must be identically $0$. 
 	 
	 \stop

 \begin{lemma}\label{explicit dbarlog} In the notation of definition \ref{ncd}, suppose that  $f:U\longrightarrow \mathbb C$ is a $\dbar\log$ compatible function which is nonzero away from $U\cap V$. Then 
 \[f= g z_{1}^{a_{1}}\dotsb z_{n}^{a_{n}}\]
 where $g:U\longrightarrow \mathbb C^{*}$ is smooth, and $a_{i}\in \mathbb N$ are locally constant.
 
\end{lemma}
 
 \pf
 
 Restrict without losing generality to the case where $\{z_{i}=0\}\subset U$ is connected for each $i$.
 Suppose that $f$ vanishes at a point $p$ where $z_{1}=0$ but $z_{i}\neq 0$ for all $i\neq 1$. We shall show that $f=hz_{1}^{n}$. There exists a continuous family of holomorphic disks intersecting $\{z_{i}=0\}$ transversely at $0$ at each point in a neighborhood of $p$ in $\{z_{i}=0\}$. Restricted to each of these holomorphic disks, $f$ is $\dbar\log$ compatible, and has an isolated zero at $0$. Let $z$ be the holomorphic coordinate on the disk. It follows from Lemma \ref{dbarlog division} and Lemma \ref{zero extension} that restricted to such a disk $f=z^{n}f'$ for some non vanishing smooth function $f'$. The number $n$ is therefore equal to the winding number of $f$ around zero in this disk. As this winding number will not change in connected families, it follows that $f$ must vanish on all of  $\{z_{1}=0\}$.   Lemma \ref{dbarlog division} then implies that $f=h z^{n}_{1}$ and that $h$ is a $\dbar\log$ compatible function which does not vanish on the set where $z_{1}= 0$ and $z_{i}\neq0$ for all $i\neq 1$.
 
 Repeating this argument, we get that $f=gz_{1}^{a_{1}}\dotsb z_{n}^{a_{n}}$ where $g$ is $\dbar\log$ compatible, and is nonzero everywhere away from where at least two of the $z_{i}$'s vanish. As this set has real codimension $4$, the winding number of $g$ around any loop must be zero, so $g$ must not vanish on any holomorphic disk with boundary outside of this set, so $g$ must be non vanishing, as required. 
 
 \stop
 
 \
 
 In the case that $M$ is holomorphic, the sheaf of holomorphic functions which are non-vanishing away from $V$ defines the log structure associated to the normal crossing divisor $V$. If the divisor is locally given by $z_{1}\dotsc z_{n}=0$, these functions are locally in the form of $gz_{1}^{a_{1}}\dotsb z_{n}^{a_{n}}$ where $g$ is holomorphic and $\mathbb C^{*}$ valued. 
 Lemma \ref{explicit dbarlog} above gives a similar property for the sheaf of $\dbar\log$ compatible functions which are non-vanishing away from $V$: any such function is locally equal to $gz_{1}^{a_{1}}\dotsb z_{n}^{a_{n}}$ where $g$ is smooth and $\mathbb C^{*}$ valued. 
 
 If we define $\mathcal N(M,V)$ to be the sheaf of $\dbar\log$ compatible functions on $M$ which are non-vanishing away from $V$, then the $C^{\infty}$ log space $(M,\mathcal N(M,V)\rightarrow \mathcal C_{M})$ is locally isomorphic to open subsets of Smooth$(X^{\dag}_{{[0,\infty)^{n}}})$. The log structure on a complex manifold with normal crossing divisors is locally isomorphic to an open subset of $X^{\dag}_{[0,\infty)^{n}}$, so in this sense, a $\dbar\log$ compatible normal crossing divisor is the $C^{\infty}$ analogue of a normal crossing divisor.

\bibliographystyle{plain}
\bibliography{ref.bib}
\end{document}